\documentclass[journal]{IEEEtran}
\usepackage[utf8]{inputenc}

\usepackage{amsmath,amsfonts,amssymb,amsthm}
\usepackage{cite}
\usepackage{array}
\usepackage{graphicx}
\usepackage{mathtools}
\usepackage{geometry}
\usepackage[dvipsnames]{xcolor}
\usepackage{soul}
\usepackage{float}
\usepackage[caption=false,font=footnotesize,labelfont=sf,textfont=sf]{subfig}
\usepackage{tikz}
\usetikzlibrary{positioning, shapes.geometric, arrows, backgrounds, decorations.pathreplacing, patterns}
\usepackage{hyperref}
\hypersetup{
    colorlinks=false,
    hidelinks,
    pdftitle={Distributed convex optimization "over-the-air" in dynamic environments.},
    pdfauthor={Navneet Agrawal, Renato L.G. Cavalcante, Masahiro Yukawa, Slawomir Stanczak},
    pdfkeywords={distributed optimization, over-the-air, OTA, decentralized federated learning}
    }
\usepackage{textcomp}

\usepackage{pgfplots}
\usepgfplotslibrary{colorbrewer,patchplots}
\pgfplotsset{compat=1.18}
\pgfplotsset{
    discard if/.style 2 args={
        x filter/.code={
            \edef\tempa{\thisrow{#1}}
            \edef\tempb{#2}
            \ifx\tempa\tempb
                
            \fi
        }
    },
    discard if not/.style 2 args={
        x filter/.code={
            \edef\tempa{\thisrow{#1}}
            \edef\tempb{#2}
            \ifx\tempa\tempb
            \else
                
            \fi
        }
    }
}

\newcommand{\Field}[1]{\mathbb{#1}}
\newcommand{\Measure}[1]{\mathbb{#1}}

\newcommand{\Set}[1]{\mathcal{#1}}

\newcommand{\Sseq}[2]{\{#1, \dots ,#2\}}

\newcommand{\Matrix}[1]{\mathbf{#1}}
\newcommand{\Vector}[1]{\pmb{#1}}
\newcommand{\Operator}[1]{\mathsf{#1}}

\newcommand{\abs}[1]{\lvert #1 \rvert}
\newcommand{\norm}[1]{\lVert#1\rVert}

\newcommand{\lev}{{\mathrm{lev}}}
\newcommand{\innerprod}[2]{\left\langle#1, #2\right\rangle}

\newcommand{\Exp}[1]{\mathbb{E}[#1]}
\newcommand{\Expp}[2]{\mathbb{E}_{#1}[#2]}
\newcommand{\cond}{\,\big| \,}

\newcommand{\probMeas}{\Measure{P}}

\newcommand{\deriv}[1]{\mathrm{d}#1}

\newcommand{\est}[1]{\hat{#1}}
\newcommand{\RE}[1]{\Re\Big[ #1 \Big]}

\newcommand{\Fix}[1]{{\mathrm{Fix}(#1)}}

\newcommand{\diagm}{\text{Diag}}

\newcommand{\mymark}[2]{\stackrel{\mathclap{\normalfont\mbox{\tiny\sffamily #1}}}{#2}}
\newcommand{\mylabel}[2]{#2\def\@currentlabel{#2}\label{#1}}

\theoremstyle{plain}
\newtheorem{theorem}{Theorem}
\newtheorem{proposition}{Proposition}
\newtheorem{lemma}{Lemma}

\theoremstyle{definition}
\newtheorem{definition}{Definition}
\newtheorem{assumption}{Assumption}

\theoremstyle{remark}
\newtheorem{remark}{Remark}

\newcommand{\Real}{\Field{R}}
\newcommand{\Complex}{\Field{C}}
\newcommand{\Natural}{\Field{N}}

\newcommand{\Realp}{\Real_+}

\newcommand{\sN}{\Set{A}}

\newcommand{\sG}{\Set{G}}
\newcommand{\sGi}{\sG_i}
\newcommand{\sGb}{\overline{\sG}_i}
\newcommand{\sE}{\Set{E}}
\newcommand{\sEi}{\sE_i}

\newcommand{\sS}{\Set{S}}

\newcommand{\sC}{\Set{C}}
\newcommand{\sCp}{\sC^\perp}
\newcommand{\sX}{\Set{X}}

\newcommand{\sD}{\Set{D}}
\newcommand{\sXki}{\sX_{k,i}}
\newcommand{\sI}{\Set{I}}

\newcommand{\sK}{\Set{K}}

\newcommand{\sM}{\Real^M}

\newcommand{\sQ}{\Set{Q}}
\newcommand{\ssT}{\Set{T}}

\newcommand{\sUp}{\Upsilon}
\newcommand{\sUki}{\Upsilon_{k,i}}
\newcommand{\sUs}{\Upsilon^\star}

\newcommand{\cT}{\Theta_{k,i}}
\newcommand{\cS}{\Theta'_{k,i}}

\newcommand{\vh}{\Vector{h}}
\newcommand{\vhki}{\vh_{k,i}}

\newcommand{\bh}{\vh^\star}

\newcommand{\vOne}{\Vector{1}}
\newcommand{\ve}{\Vector{e}}
\newcommand{\vn}{\Vector{n}}
\newcommand{\vv}{\Vector{v}}

\newcommand{\vni}{\vn_i}
\newcommand{\vti}{\tilde{\vn}_i}
\newcommand{\vy}{\Vector{y}}

\newcommand{\vnn}{\acute{\vn}}

\newcommand{\vx}{\Vector{x}}
\newcommand{\vxs}{\vx^\star}

\newcommand{\vu}{\Vector{u}}
\newcommand{\vz}{\Vector{z}}
\newcommand{\vzki}{\vz_{k,i}}
\newcommand{\vzero}{\Vector{0}}

\newcommand{\gP}{\Vector{\psi}}
\newcommand{\gPi}{\gP_i}
\newcommand{\gZ}{\Vector{z}}
\newcommand{\gZi}{\gZ_i}
\newcommand{\gZki}{\gZ_{k,i}}

\newcommand{\gPs}{\gP^\star}
\newcommand{\gH}{\Vector{\phi}}
\newcommand{\gHi}{\gH_i}

\newcommand{\gA}{\Vector{\alpha}}

\newcommand{\gAki}{\gA_{k,i}}

\newcommand{\gL}{\Vector{\lambda}}
\newcommand{\gLi}{\gL_i}
\newcommand{\gLki}{\gL_{k,i}}

\newcommand{\mX}{\Matrix{X}}
\newcommand{\mJ}{\Matrix{J}}
\newcommand{\mZ}{\Matrix{Z}}

\newcommand{\mI}{\Matrix{I}}
\newcommand{\mL}{\Matrix{L}}
\newcommand{\mLi}{\mL_i}
\newcommand{\mLbi}{\overline{\mL}_i}

\newcommand{\mP}{\Matrix{P}}
\newcommand{\mPi}{\mP_i}
\newcommand{\mPbb}{\overline{\mP}}
\newcommand{\mPb}{\mPbb_i}
\newcommand{\mG}{\Matrix{G}}
\newcommand{\mGi}{\mG_i}
\newcommand{\mGb}{\overline{\mG}_i}
\newcommand{\mZero}{\Matrix{0}}

\newcommand{\mW}{\Matrix{W}}
\newcommand{\mWi}{\mW_i}

\newcommand{\opT}{\Operator{T}}
\newcommand{\opTi}{\opT_i}
\newcommand{\opTki}{\opT_{k, i}}

\newcommand{\opP}{\Operator{P}}

\newcommand{\opL}{\Operator{K}}
\newcommand{\opLi}{\opL_i}

\newcommand{\tZ}{\tilde{\rnu}}

\newcommand{\dmin}{\delta_{\min}}
\newcommand{\dmax}{\delta_{\max}}

\newcommand{\tW}{\tilde{w}}
\newcommand{\sif}{\sum_{i=0}^\infty}
\newcommand{\iN}{{i\in\Natural}}
\newcommand{\nN}{{n\in\Natural}}

\newcommand{\bB}{{b\in\Sseq{1}{2B}}}
\newcommand{\kN}{{k\in\sN}}
\newcommand{\oinO}{{\omega\in\Omega}}

\newcommand{\sigmaF}{\mathcal{F}}
\newcommand{\pspace}{(\Omega, \sigmaF, \mathbb{P})}

\newcommand{\VS}{\Real^{MN}}
\newcommand{\VM}{\Real^{MN \times MN}}

\newcommand{\rmuki}{\mu_{k,i}}
\newcommand{\pPi}{\mathbb{P}_i}

\newcommand{\oP}{\check{\mP}}
\newcommand{\oW}{\check{\mW}}

\newcommand{\ovn}{\check{\vn}}
\newcommand{\oPi}{\oP_i}

\newcommand{\oWi}{\oW_i}

\newcommand{\ovni}{\ovn_i}
\newcommand{\ogL}{\check{\gL}}
\newcommand{\ogP}{\check{\gP}}

\newcommand{\rnu}{\nu}
\newcommand{\rz}{\rnu_{jr}}
\newcommand{\rni}{\rnu_{jr,i}}

\newcommand{\eNoise}{m_{\tW}}

\newcommand{\hf}{\hat{f}}
\newcommand{\Ker}{\Operator{k}}
\newcommand{\ey}{\est{y}}
\newcommand{\sQki}{\sQ_{k,i}}

\newcommand{\ow}{\Vector{\varpi}}

\newcommand{\Epp}[1]{\Exp{#1\, \big|\, \gPi}}
\newcommand{\Epl}[1]{\Exp{#1\, \big|\, \gP_l}}

\newcommand{\fs}{f^{\star}}
\newcommand{\fsp}{\fs(\gP_{i+1})}
\newcommand{\fse}{\fs(\gPi)}
\newcommand{\fq}{\Operator{Q}}

\newcommand{\pv}{\vartheta}
\newcommand{\pvv}{\Vector{\pv}}
\newcommand{\pb}{\theta}

\newlength{\bulletgap}
\begin{document}

\title{Distributed Convex Optimization ``Over-the-Air'' in Dynamic Environments}
\author{Navneet~Agrawal,~Renato~L.G.~Cavalcante,~\IEEEmembership{Member,~IEEE,} Masahiro~Yukawa,~and~S{\l}awomir~Sta\'nczak,~\IEEEmembership{Senior Member,~IEEE}%
        \thanks{Navneet~Agrawal is with Network Information Theory Group at Technische Universit\"at Berlin, Germany (navneet.agrawal@tu-berlin.de), and his research is funded by the German Research Foundation (DFG) within their priority program SPP1914 ``Cyber-Physical Networking''.}%
        \thanks{Renato Lu\'is Garrido Cavalcante and S{\l}awomir Sta\'nczak are with Department of Wireless Communications and Networks at Fraunhofer Heinrich Hertz Institute, Germany, and they acknowledge the support of joint project 6G-RIC with project identification numbers 16KISK020K and 16KISK030.}%
        \thanks{Masahiro Yukawa is with Department of Electronics and Electrical Engineering, Keio University, Japan, supported by JST SICORP under Grant JPMJSC20C6, Japan.}%
        }

\maketitle

    \begin{abstract}

This paper presents a decentralized algorithm for solving distributed convex optimization problems in dynamic networks with time-varying objectives.
The unique feature of the algorithm lies in its ability to accommodate a wide range of communication systems, including previously unsupported ones, by abstractly modeling the information exchange in the network.
Specifically, it supports a novel communication protocol based on the ``over-the-air'' function computation (OTA-C) technology, that is designed for an efficient and truly decentralized implementation of the consensus step of the algorithm.
Unlike existing OTA-C protocols, the proposed protocol does not require the knowledge of network graph structure or channel state information, making it particularly suitable for decentralized implementation over ultra-dense wireless networks with time-varying topologies and fading channels.
Furthermore, the proposed algorithm synergizes with the ``superiorization'' methodology, allowing the development of new distributed algorithms with enhanced performance for the intended applications.
The theoretical analysis establishes sufficient conditions for almost sure convergence of the algorithm to a common time-invariant solution for all agents, assuming such a solution exists.
Our algorithm is applied to a real-world distributed random field estimation problem, showcasing its efficacy in terms of convergence speed, scalability, and spectral efficiency. 
Furthermore, we present a superiorized version of our algorithm that achieves faster convergence with significantly reduced energy consumption compared to the unsuperiorized algorithm.
    \end{abstract}

    \section{Introduction}\label{sec:intro}
\IEEEPARstart{C}{ommunication} and computation are both fundamental in realizing multiagent systems where interconnected devices cooperate to solve inference, control, or learning tasks~%
\cite{nedic2009distributed,srivastava2011distributed,cavalcante2009adaptive,cavalcante_dynamic_2013,nedic2018distributed}.
In such systems, the data associated with the tasks are distributed among the devices, which take advantage of their aggregated computational power by exchanging information, without resorting to a unique (centralized) provider.
Often, communication is the main bottleneck for implementing distributed algorithms in these systems.
Specifically, the densification of the network, coupled with the scarcity of communication resources, poses severe scalability issues, especially in time-varying environments.
Hence, the development of novel distributed algorithms that support robust, scalable, and energy-efficient communication is of fundamental importance for next generation wireless multiagent systems \cite{notarstefano2019distributed}.

\subsection{Distributed systems with time-varying objectives}

In this study, the objective of the proposed algorithm is to minimize the sum of local costs under local constraints, each of them known only to an individual agent~%
\cite{nedic2009distributed,srivastava2011distributed,cavalcante2009adaptive,cavalcante_dynamic_2013,nedic2018distributed}.
In addition, we address problems in which the optimization objectives are possibly nonsmooth and time-varying, under the assumption that a time-invariant solution exists. 
Problems of this type are common in wireless sensor networks or autonomous systems, where each sensor continuously obtains a potentially infinite number of measurements of a common quantity of interest at regular intervals.
Concrete applications include, among others, decentralized source localization~\cite{cavalcante_dist_2011}, signal detection \cite{cavalcante_dynamic_2013}, and federated learning~\cite{konevcny2016federated,sahin2022distributed,michelusi2022non}.

Solvers addressing these distributed optimization problems typically apply a variation of the following iterative two-step approach~%
\cite{nedic2009distributed,boyd2011distributed,duchi2011dual,nedic2018distributed,yang2019survey}:
the \emph{local optimization} step, where agents update their estimate of a global optimizer locally,
followed by the \emph{consensus} step, where agents seek agreement by fusing their local estimates with their neighbors.
In the local optimization step, agents only utilize information available locally to them at the time, such as measurements and/or previous estimates of the quantity of interest.
Typically, it is implemented using some gradient-based approach \cite{boyd2011distributed,duchi2011dual}, and in this paper, we employ the adaptive projected subgradient method (APSM)\cite{yamada2005adaptive,Slavakis_2006}.
The APSM provides an adaptive solution to a sequence of problems over an infinite time horizon, and it has been successfully applied to problems in both centralized and distributed settings with nonsmooth and time-varying objectives, as seen in previous studies~\cite{cavalcante2009adaptive,cavalcante_dist_2011,cavalcante_dynamic_2013,fink2022superiorized}.

In the consensus step, agents employ a consensus protocol that relies on information from their neighboring agents.
This information is obtained through a communication mechanism over the wireless network. 
Most existing algorithms utilize a communication mechanism based on the channel separation principle, such as gossip or broadcast communication.
The amount of wireless resources required in such a communication grows proportionally with the increase in number of agents sharing the channel.
Furthermore, channel separation-based communication is known to be inefficient when the objective is to evaluate a function of the data distributed across the network \cite{nazer2007computation,goldenbaum2013robust}.
This inefficiency can be addressed with the \emph{over-the-air function computation} (OTA-C) technology \cite{goldenbaum2013robust,frey_ota_2019}, which has shown excellent performance in various consensus-based applications \cite{abari2016over,agrawal2019scalable,molinari2021max,molinari2022over}.
In OTA-C, the processes of function computation and data transmission are merged in the communication system by exploiting the additive structure of the wireless multiple access channels, and hence, the wireless resources required by the OTA-C does not increase with the growing network density.
However, existing OTA-C protocols \cite{goldenbaum2013robust,abari2016over,frey_ota_2019,agrawal2019scalable,molinari2021max,molinari2022over} are designed for centralized network architectures, where a central node is required to obtain, for example, channel state information (CSI) of individual links.
Hence, existing OTA-C protocols can be challenging to implement in networks that are ultra-dense, decentralized, and time-varying.

\subsection{Summary of contributions}
Building upon \cite{cavalcante_dynamic_2013}, we propose a novel distributed algorithm for solving convex optimization problems with time-varying objectives over decentralized wireless networks.
The proposed algorithm supports a wide range of communication mechanisms by modeling information exchange for the consensus step in an abstract fashion.
As a result, the algorithm encompasses the communication employed in previous studies \cite{nedic2009distributed,boyd2011distributed,cavalcante_dynamic_2013} as a particular case.
Moreover, it includes the novel OTA-C protocol proposed in this study, as detailed in Section \ref{sec:comm}, which was previously unsupported.

One of the key contributions of this paper is development of a novel OTA-C protocol, which addresses several important challenges in implementing the consensus step over wireless networks. 
The protocol can handle fading, noisy, and time-varying wireless channels.
Moreover, it stands out as a \emph{truly} decentralized solution for the following reasons:
(1) it is robust against node failures, 
(2) it does not rely on prior coordination among agents, and 
(3) its communication overhead remains independent of the network size.
The proposed protocol stands out from existing OTA-C protocols \cite{goldenbaum2013robust,abari2016over,frey_ota_2019} due to its unique design that enables operation without requiring information such as CSI that could potentially impede a truly decentralized implementation.

Our main theoretical results (Theorem \ref{thm:main}) establish sufficient conditions for the asymptotic convergence of the proposed algorithm to a common feasible solution for all agents, assuming that at least one such a solution exists.
We provide \emph{almost sure} convergence guarantees for our algorithm, which stands apart from the guarantees of convergence in the mean square sense, provided by most existing studies \cite{nedic2009distributed,duchi2011dual,michelusi2022non} (definitions of convergence in almost sure and mean square sense can be found in Section \ref{sec:prelims}).
In addition, our theoretical results hold under weaker conditions (see Assumption \ref{ass:random_ass} and Remark \ref{rem:comm_ass}) than those imposed in the previous studies \cite{nedic2009distributed,cavalcante_dynamic_2013} tackling similar problems.

The proposed algorithm is proven to be resilient to random perturbations, meaning that the algorithm is guaranteed to converge to a solution of the original problem, even if random disturbances (in the sense of Definition \ref{def:perturb}) are present in its iterates.
Furthermore, leveraging this resilience property in Section \ref{sec:sim}, we design a ``superiorized'' version of our algorithm by adding carefully designed perturbations to its iterates that enhances the performance for the intended application compared to its unperturbed version \cite{censor2010perturbation,fink2022superiorized}.
To the best of our knowledge, this algorithm is the first to integrate the \emph{superiorization} methodology \cite{censor2010perturbation,fink2022superiorized} in a truly decentralized context.

We apply the proposed algorithm to a real-world problem of distributed random field estimation in Section \ref{sec:sim}, and demonstrate its efficacy in terms of convergence speed, scalability, and spectral efficiency.
Moreover, inspired by \cite{candes2008enhancing,pollakis2012base,fink2022superiorized}, we present a superiorized version of our algorithm that is designed to be faster and more energy-efficient.
Specifically, the new algorithm leads to a significant reduction in the energy consumption in communication (over 80\% in our experiments), while maintaining only a small impact on the algorithm's performance.\footnotemark

\footnotetext{Part of this study has been recently published in a conference \cite{agrawal_icassp23}.
However, unlike \cite{agrawal_icassp23}, in this manuscript, we establish the bounded perturbation resilience of the proposed algorithm, and, in Section \ref{sec:sim}, this property is used to design a new energy-efficient version of the algorithm.
Moreover, this manuscript includes proofs of theorems that are more general than those in the conference manuscript \cite{agrawal_icassp23}.}

\subsection{Related existing studies}
Recently, OTA-C based distributed optimization algorithms that do not necessarily require a centralized architecture have been proposed in \cite{xing2021federated,amiri2021blind,sahin2022distributed,michelusi2022non}.
However, these are not truly decentralized as they require some level of prior coordination (e.g.~graph coloring in \cite{xing2021federated}), or additional overheads (e.g.~beamforming in \cite{amiri2021blind}), or their algorithm is only suitable for a specific optimization scenario (e.g.~signSGD with majority voting consensus in \cite{sahin2022distributed}).
In \cite{michelusi2022non}, the author proposes a truly decentralized algorithm, but only considers problems with time-invariant and smooth objectives.
Moreover, the OTA-C protocol in \cite{michelusi2022non} exacerbates the effective noise in the system, and the sufficient conditions for convergence of the algorithm require strong assumptions on the communication system.
In contrast, our results are established for weaker assumptions on the channel (see Assumption \ref{ass:WMAC}). 
Hence, our OTA-C protocol can be deployed on very general and realistic wireless communication scenarios.

\vspace{1em}

The paper is organized as follows:
basic notations, definitions, and preliminary results are presented in Section \ref{sec:prelims}.
Section \ref{sec:prob} introduces the system and the class of distributed optimization problems considered in this paper.
In Section \ref{sec:alg}, we describe our proposed algorithm, and establish sufficient conditions for asymptotic convergence of the scheme to a feasible solution.
The novel OTA-C protocol is presented in Section \ref{sec:comm}, where we also investigate conditions for its applicability to our proposed scheme.
An application of the proposed scheme is presented in Section \ref{sec:sim} for the problem of distributed supervised random field estimation.

    \section{Notations and Preliminaries}
    \label{sec:prelims}



We use symbols $\Real$, $\Realp$ and $\Natural$ for real, nonnegative real, and natural (including zero) numbers, respectively.
Scalars, vectors, and matrices are typified with normal, bold, and capital-bold typeface characters.
The $(p,q)$ component of a matrix $\mX$ is denoted by $[\mX]_{(p,q)}$.
A sequence is denoted by $(x_i)_{i\in\sI}$ for some index set $\sI$, or by simply $(x_i)$ when the index set is clear from the context.
A diagonal matrix with diagonal elements corresponding to the elements of the vector $\vx$ is denoted as $\diagm(\vx)$.

In this paper, any real-valued $M$-dimensional vector is an element of the Hilbert space $(\Real^M, \innerprod{\cdot}{\cdot})$ with inner-product $(\forall \vv,\vy \in \Real^M)$ $\innerprod{\vv}{\vy} := \vv^T\vy$ and induced norm $\norm{\vv} := \sqrt{\vv^T \vv}$.
Given a matrix $\mX\in\Real^{M\times N}$, its spectral norm is $\norm{\mX}_2 := \max_{\norm{\vy} \neq 0} \norm{\mX\vy} / \norm{\vy}$,
and its spectral radius is $\rho(\mX) := \max_i\{\abs{\lambda_i} \ \big|\ \lambda_i \text{ is an eigenvalue of } \mX\}$.
Matrices satisfying $\rho(\mX) = \norm{\mX}_2$ are called the \emph{radial} matrices \cite{goldberg1974radial}. 
In particular, all symmetric ($\mX^T = \mX$) and normal ($\mX\mX^T = \mX^T\mX$) matrices are radial.

By $\pspace$ we denote a probability space, where $\Omega$ is an arbitrary sample set, $\mathcal{F}$ is a $\sigma$ algebra of subsets of $\Omega$, and $\probMeas$ is a probability measure.
As common in literature \cite{shapiro2021lectures}, we define a \emph{random variable} as a measurable map $x : \Omega \to \Real$.
A random quantity $x$ is sometimes viewed as a map, e.g.~$\Omega\to\Real$, and sometimes as its particular realization, i.e.~$x := x(\omega) \in \Real^M$.
The meaning will be clear from the context.

A directed graph of networked agents is represented by $\sG=(\sN, \sE)$, where $\sN := \{1,\dots,\abs{\sN}\}$ is the index set of agents, and $\sE \subseteq \sN\times\sN$ contains ordered pairs of connected agents such that if $(k,l)\in\sE$, then agent $k$ can transmit to agent $l$.
Here, $\abs{\sN}$ denotes cardinality of the set $\sN\subseteq\Natural$.
A \emph{weighted} graph has each edge $(k,l)\in\sN\times\sN$ endowed with a weight $\alpha_{k,l}\in\Real$ such that $\alpha_{k,l} \neq 0$ if $(k,l)\in\sE$, and zero otherwise.
A directed graph becomes \emph{undirected} if for every $(k,l)\in\sE$, we also have $(l,k)\in\sE$, and $\alpha_{k,l} = \alpha_{l,k}$.
An undirected graph is called \emph{connected} if, for any $(i,j)\in\sN\times\sN$, there is a sequence of edges in the set $\sE$ that make a path between $i$ and $j$.
A \emph{random graph} is a weighted graph with edge weights $(\alpha_{k,l})$ being random variables.
The \emph{expected graph} of a random graph $\sG$, denoted by $\overline{\sG}$, is the one with edge weights as the expected weights $(\Exp{\alpha_{k,l}})$, given that the expected values exists.

A set $\sC \subseteq \Real^M$ is said to be \emph{convex} if $\lambda \vv + (1 - \lambda) \vy \in \sC$ for every $\vv\in\sC$, $\vy\in\sC$ and $\lambda \in (0, 1)$.
A \emph{closed} set is the one that also contains all its boundary points.
Let $\sX$ be a nonempty closed convex set in $\Real^M$.
Then, by $\opP_\sX : \Real^M \to \sX$, we define an \emph{orthogonal projection mapping} \cite{Bauschke_2017} that maps any point $\vv\in\Real^M$ to a unique point $\opP_\sX(\vv)\in\sX$ such that $\norm{\vv - \opP_\sX(\vv)} = \min_{\vy\in\sX} \norm{\vv - \vy}$.
Moreover, a mapping $\opT : \sX \to \Real^M$, with nonempty fixed-point set $\Fix{\opT} := \{ x\in\sX \mid \opT(x) = x\}$, is called \emph{quasi-nonexpansive} if $(\forall \vx\in\sX)$ $(\forall\vy\in\Fix{\opT})$ $\norm{\opT(\vx) - \vy} \leq \norm{\vx - \vy}$.
For such a mapping $\opT$, the fixed-point set $\Fix{\opT}$ is also closed and convex \cite[Corollary 4.24]{Bauschke_2017}.
The \emph{subdifferential} of a function $\Theta:\Real^M\to\Realp$ at a point $\vy\in\Real^M$ is a closed convex set $\partial\Theta(\vy) := \{\vh \in \Real^M \mid \quad \Theta(\vy) + \innerprod{\vx - \vy}{\vh} \leq \Theta(\vx), \forall \vx\in\Real^M\}$.
A \emph{subgradient} of $\Theta$ at $\vy$, denoted by $\Theta'(\vy)$, is a particular selection of element from $\partial\Theta(\vy)$.
The \emph{$c$-sublevel set} of a convex function $\Theta$, defined as $\lev_{\leq c} \Theta := \{\vh\in\Real^M \mid \Theta(\vh) \leq c\}$, is a convex set for every $c\in\Realp$ \cite{Bauschke_2017}.

Let $(X_n)_\nN$ be a sequence of random numbers taking values in $\sX\subset \Real^M$.
Then, $(X_n)_\nN$ is said to converge to $X\in\sX$ (a) in the \emph{mean square sense} if $\lim_{n\to\infty}\Exp{\norm{X_n - X}^2} = 0$, and 
(b) \emph{almost surely} (a.s) if $\probMeas(\lim_{n\to\infty} X_n = X) = 1$.
Moreover, $(X_n)_\nN$ is said to be \emph{bounded a.s.} if for \emph{almost every} (a.e.) $\oinO$ (i.e.~excluding points in $\Omega$ with a measure zero) there exists $b(\omega)\in\Realp$ such that $X_n(\omega) \leq b(\omega)$ for all $\nN$.

\begin{definition}[Sequence of bounded perturbations] \label{def:perturb}
    Given a sequence of nonnegative scalars $(\alpha_n)_{\nN}$, and a sequence of random vectors $(\vy_n)_{\nN}$, the sequence $(\alpha_n \vy_n)_\nN$ is called a \emph{sequence of bounded perturbations} if the series $\sum_\nN \alpha_n$ is convergent, and the sequence $(\vy_n)_\nN$ is bounded a.s.
\end{definition}

    \section{System and Problem Description}
    \label{sec:prob}


Let $\iN$ denote the time-index, and define $\sGi := (\sN, \sEi)$ as the graph formed by $N:=\abs{\sN}$ agents at time $i$.
We consider the following problem $\pPi$ at any time $i$:
\begin{equation}
    \label{eq:optprob}
    \begin{aligned} 
        \underset{\gP\in\sX_i}{\text{minimize}} \ \Theta_i(\gP) := \sum_{k\in\sN} \cT(\vh_k), \quad
        \text{s.t.}\ \gP\in\sC,
    \end{aligned}
\end{equation}
where, given $(\forall \kN)\ \vh_k\in\Real^M$, the variable $\gP\in\Real^{MN}$ takes the form $\gP:=[\vh_1^T, \dots, \vh_N^T]^T$, the local cost function of any agent $\kN$ is defined as $\cT:\Real^M\to\Realp$, and, given $\sXki \subseteq \Real^M$ as a local constraint set of any agent $\kN$, the set $\sX_i\subseteq \Real^{MN}$ is defined as $\sX_i := \sX_{1,i} \times \dots \times \sX_{N, i}$.
The \emph{consensus subspace} $\sC$, defined in the following, plays an important role in this study.

\begin{definition}[Consensus subspace $\sC$]
    The \emph{consensus subspace} $\sC \subset \VS$ is defined as
    \begin{equation}
    \label{eq:consensus_subspace}
    \begin{aligned}
        \sC := \{[\vh_{1}^T, \dots, \vh_{N}^T]^T \mid \ &\vh_1 = \dots = \vh_N, \\
            & (\kN)\ \vh_k\in\Real^M\}.
    \end{aligned}
    \end{equation}
    The orthogonal projection mapping onto $\sC$ is given by $\opP_{\sC}(\vx) := \mJ \vx$, where $\mJ := [\vu_1, \dots, \vu_M] [\vu_1, \dots, \vu_M]^T \in \VM$.
    The vector $\vu_j$, for $j=1,\dots,M$, is defined as $\vu_j := (\vOne_N \otimes \ve_j)/\sqrt{N} \in \VS$, where $\vOne_N\in\Real^N$ is the vector of ones, $\ve_j\in\Real^M$ is the standard $j$th basis vector, and $\otimes$ denotes the Kronecker product.
\end{definition}

%
\begin{assumption}
    \label{ass:objective}
        (i) For all $\kN$ and $\iN$, $\cT$ is a continuous convex function, and $\sXki$ is a closed convex set. 
            Both $\cT$ and $\sXki$ are assumed to be known to agent $k$ at time $i$.
            Moreover, the set $\sXki$ is characterized as the fixed-point set of a quasi-nonexpansive operator [defined later in \eqref{eq:local_scheme}], and contains minimizers of $\cT$ that attain the value of zero on $\cT$.
            Formally, $(\forall \kN)(\forall \iN)$
            \begin{align}\label{eq:solution_set}
                \sUki := \{\vh\in\sXki \mid \cT(\vh) = 0\} \neq \emptyset.
            \end{align}

        (ii) Every set in the sequence $(\sUki)_{\kN,\iN}$ contains at least one point that is common to all.
            More precisely, $\sUs := \bigcap_{\iN} \bigcap_{k\in\sN} \sUki \neq \emptyset$.

\end{assumption}

\begin{remark} \label{rem:prob_form}
    Many multi-agent applications require agents to collaboratively find a common solution from the (nonempty) set of feasible solutions to all problems in the sequence $(\pPi)_\iN$. 
    Assumption \ref{ass:objective} is standard in such applications \cite{yamada2005adaptive,cavalcante2009adaptive}.
    Examples of concrete applications of this type include collaborative source localization \cite{cavalcante_dist_2011}, distributed detection \cite{cavalcante_dynamic_2013}, and decentralized federated learning \cite{konevcny2016federated}.
    Assumption \ref{ass:objective}(i) may seem to be rather restrictive due to \eqref{eq:solution_set}, but it is often easily satisfied for many set-theoretic problems with available knowledge about the system.
    Nonetheless, we show in Section \ref{sec:sim} that our proposed method works well even if Assumption \ref{ass:objective} is not fully satisfied because of, for instance, noise in the measurements.
\end{remark}

Ideally, our objective is to design a distributed algorithm that allows all agents $k\in\sN$ to find the same vector in the set $\sUs$.
As discussed in \cite{cavalcante_dynamic_2013}, finding a common solution to all problems in the sequence $(\pPi)_\iN$ is challenging because of, for instance, system causality (the $i$th measurement is only available after time $i$ has elapsed)
and memory limitations (previous information must be discarded to free memory for new information).
Dynamical systems of this type are typically approached by solving each problem in the sequence $(\pPi)_\iN$ independently by using local communication mechanisms to achieve consensus.
In large distributed systems, where agents can only communicate with a few neighbors, this approach may take too many iterations to solve any single problem instance.
Moreover, independently solving the problems ignores the rich temporal structure in the sequence of information available over time, and it does not necessarily lead to a good estimate of the desired parameter.
Indeed, the set of solutions to any individual problem $\pPi$ may contain points that are arbitrarily far from the desired true parameter, e.g.~when the solution set of $\pPi$ is a linear variety in $\Real^M$.
Hence, without taking into account all problems available from the past, each individual problem $\pPi$ provides little information about the \emph{estimandum}.
Therefore, in order to obtain a good estimate, distributed and adaptive algorithms aim to find a point that simultaneously solves as many problems $\pPi$ in \eqref{eq:optprob} as possible, under Assumption \ref{ass:objective}(ii), which guarantees that a common solution to all problems $(\pPi)_\iN$ exists.
Building upon \cite{cavalcante_dynamic_2013}, our objective in Section \ref{sec:alg} is to make agents agree on a time-invariant solution to all \emph{but finitely many} problems in $(\pPi)_\iN$.
More precisely, agents must agree on a solution to the following feasibility problem with infinitely many constraints:
\begin{align}
    \label{eq:relax_set}
    \text{Find}\ \bh \in \sUp := \overline{\liminf_{i\to\infty} \sUp_i} = \overline{\bigcup_{\nN}\bigcap_{i\geq n} \sUp_i} \supset \sUs,
\end{align}
where the overline denotes closure of the set, $\sUp_i := \bigcap_{k\in\sN} \sUki$, and $\sUki$ is defined in \eqref{eq:solution_set}.
    
    \section{Proposed Algorithm}
    \label{sec:alg}

In this section, we propose a distributed algorithm to find a solution of the feasibility problem \eqref{eq:relax_set}.
As discussed in Section \ref{sec:intro}, the proposed algorithm follows an iterative two-step approach consisting of a \emph{local optimization step}, followed by a \emph{consensus step}~\cite{nedic2009distributed,cavalcante2009adaptive,cavalcante_dynamic_2013}.
In the local optimization step, we adopt the well-known APSM.
Moreover, to establish the \emph{bounded perturbation resilience} of the algorithm, we introduce a sequence of random perturbations (in the sense of Definition \ref{def:perturb}) to the iterates.
In an application of the proposed algorithm in Section \ref{sec:sim}, we leverage this property by adding carefully designed perturbations to obtain sparse solutions within the feasible set, resulting in enhanced performance, and a significant reduction in energy consumption during communication.

In the consensus step of our algorithm, we employ a first-order consensus protocol, given in \eqref{eq:local_scheme_consensus} or \eqref{eq:scheme_consensus}, where information exchange among all agents in the network is described by an abstract mathematical model \eqref{eq:comm_model} that captures the underlying communication protocol.
The proposed abstraction covers a large class of wireless communication protocols, many of which are unsupported by the models proposed in the previous studies \cite{nedic2009distributed,cavalcante2009adaptive,cavalcante_dynamic_2013}.
In particular, the proposed algorithm can be successfully applied to systems employing the OTA-C protocol proposed in Section \ref{sec:comm}, which is not the case with the model proposed in \cite{cavalcante_dynamic_2013}.
In the following, we use the notation $x = x(\omega)$ to emphasize that $x$ is a random quantity.

Motivated by \cite{cavalcante_dynamic_2013} and references therein, the proposed algorithm generates a stochastic sequence $(\vhki)_\iN$ via: $(\forall \kN)(\forall \iN)(\oinO)$
\begin{subequations}
    \label{eq:local_scheme}
    \begin{align}
        &\gLki(\omega) = \opTki \big( \vhki(\omega) - \gAki(\omega) + \zeta_i\ \gZki(\omega) \big), \label{eq:local_scheme_local_opt} \\
        &\vh_{k, i+1}(\omega) = (1 - \beta_i) \gLki(\omega) + \beta_i \opL_{k,i}\big(\omega, \gLi(\omega)\big), \label{eq:local_scheme_consensus}
    \end{align}
\end{subequations}
where the starting point $\vh_{k,0}\in\Real^M$ can be arbitrary.
The mapping $\opTki$ is defined as $\opTki:\Real^M \to \Real^M$;
$(\zeta_i \gZki(\omega))_\iN$ is a sequence of bounded perturbations in the sense of Definition \ref{def:perturb}; the scalar $\rmuki \in (0,2)$ is a design parameter;
the vector $\gAki(\omega)$ is defined as
    \begin{align}\label{eq:alpha}
        \gAki(\omega) := \rmuki\ \frac{\cT(\vhki(\omega))}{\norm{\cS(\vhki(\omega))}^2}\  \cS(\vhki(\omega)),
    \end{align}
if $\norm{\cS(\vhki(\omega))}^2 \neq 0$, otherwise $\gAki(\omega) = 0$,
where $\cS(\vhki(\omega)) \in\partial \cT(\vhki(\omega))$; 
the vector $\gLi(\omega)$ takes the form $[\gL_{1,i}^T(\omega), \dots, \gL_{N,i}^T(\omega)]^T$;
and the random mapping $\opL_{k,i}:\Omega\times\Real^{MN} \to \Real^M$ models the information obtained by agent $k$ from other agents over the wireless network.
Equation \eqref{eq:local_scheme_local_opt} is the local optimization step, and equation \eqref{eq:local_scheme_consensus} is the consensus step.
Assumption \ref{ass:local_scheme} below characterizes the sequence $(\beta_i)_\iN$, the mapping $\opTki$, and the random operator $\opL_{k,i}$.

\begin{assumption} 
    \label{ass:local_scheme}
    For all agents $\kN$ using the scheme in \eqref{eq:local_scheme} at any time $\iN$, the following properties hold:

        (i) The mapping $\opTki$ is quasi-nonexpansive with $\Fix{\opTki} = \sXki$, where $\sXki$ is known to agent $k$ at time $i$ (see Assumption \ref{ass:objective}(i)).
        
        
            
        (ii) The sequence $(\beta_i)_\iN$ is a design parameter such that $(\forall \iN)$ $\beta_i \in (0, 1]$, the series $\sum_\iN \beta_i$ diverges to infinity, and the series $\sum_\iN \beta^2_i$ converges.
        
        (iii) The random operator $\opL_{k,i}$ has the form:
            \begin{align} \label{eq:comm_local}
                \opL_{k,i}(\omega, \gL) := (\mP_{k,i}(\omega) + \mW_{k,i}(\omega))\ \gL + \vn_{k,i}(\omega)
            \end{align}
            where $\mP_{k,i}$, $\mW_{k,i}$, and $\vn_{k,i}$ take values in $\Real^{M\times MN}$, $\Real^{M\times MN}$, and $\Real^M$, respectively.
            

\end{assumption}

We omit $\omega$ in the following for clarity. 
Henceforth, define $\gLi$, $\gPi$, $\gHi$, and $\gZi$ as vectors in $\Real^{MN}$ formed by stacking the vectors $\gLki$, $\vhki$, $\gAki$, and $\gZki$, respectively, for all $\kN$.
Similarly, define mappings $\opTi$ and $\opLi$ such that, for any $\gP:=[\vh_1^T,\dots,\vh_N^T]^T$ and $\gL\in\Real^{MN}$, the vectors $\opTi(\gP)\in\Real^{MN}$ and $\opLi(\gL)\in\Real^{MN}$ are formed by stacking $\opTki(\vh_k)\in\Real^{M}$ and $\opL_{k,i}(\gL)\in\Real^{M}$, respectively, for all $\kN$.
For the subsequent treatment and analysis, we represent the scheme \eqref{eq:local_scheme} as: $(\forall \iN)$
\begin{subequations}
    \label{eq:scheme}
    \begin{align}
        \gLi &= \opTi \big( \gPi - \gHi + \zeta_i\ \gZ_i \big), \label{eq:scheme_local_opt} \\
        \gP_{i+1} &= (1 - \beta_i)\ \gLi + \beta_i\ \opLi\big(\gLi\big), \label{eq:scheme_consensus}
    \end{align}
\end{subequations}
where the mapping $\opLi$ is given by:
\begin{align}
    \label{eq:comm_model}
    \opLi(\gL) := \big(\mPi + \mWi \big)\, \gL + \vni.
\end{align}
The quantities $\mPi$, $\mWi$, and $\vni$ are formed by stacking the rows of $\mP_{k,i}$, $\mW_{k,i}$, and $\vn_{k,i}$, respectively, for all agents $\kN$.

\begin{remark}
    Note that the matrix $\mPi$ in \eqref{eq:comm_model} takes the role of a random consensus matrix, defined below in Definition \ref{def:consensus_matrix}.
    The motivation for the model in \eqref{eq:comm_model} comes from the proposed OTA-C protocol, described later in Section \ref{sec:comm}, where the information exchange model [see \eqref{eq:ota_protocol}] is a particular instance of \eqref{eq:comm_model}.
    In that model, components of $\mWi$ and $\vni$ correspond to noise at the receiver, and components of $\mPi$ to fading channels.
    It is also worth noting that such communication systems were not supported by the models used in the previous studies \cite{nedic2009distributed,boyd2011distributed,cavalcante_dynamic_2013}.
\end{remark}

Assumption \ref{ass:random_ass} below, along with Definition \ref{def:consensus_matrix}, formalizes conditions imposed on the model in \eqref{eq:comm_model}.
\begin{definition}[$(\epsilon, \delta)$-random consensus matrix]
    \label{def:consensus_matrix}
    \hfill
    
    For given $0 < \epsilon \leq 1$ and $0 \leq \delta < \infty$, a random matrix $\mP$ is said to be an \emph{$(\epsilon, \delta)$-random consensus matrix} if it satisfies the following properties:

        (i) $\mP \vv = \vv$ for every $\vv \in \sC$ a.s.
        
        (ii)  $\mPbb := \Exp{\mP}$ is a \emph{nonnegative matrix}.
        
        (iii) 
            With $\sCp := \{\vx \in \VS \mid (\forall \vy \in \sC)\ \vx^T\vy = 0\}$,
            \begin{align} \label{eq:P_conn}
                \min_{\norm{\vx}=1,\, \vx \in \sCp} \vx^T \mPbb \vx \leq 1 - \epsilon.
            \end{align}
            
        (iv) $\mPbb$ is a \emph{radial} matrix \cite{goldberg1974radial}; i.e.~$\rho(\mPbb) = \norm{\mPbb}_2$.
        
        (v) $1 \leq \norm{\Exp{\mP^T \mP}}_2 \leq 1 + \delta$.

\end{definition}

\begin{assumption} 
    \label{ass:random_ass}
    For every $\iN$, the following holds:
    
        (i) $\mPi$ is an \emph{$(\epsilon_0, \delta_0)$-random consensus matrix} in the sense of Definition \ref{def:consensus_matrix}.
    
        (ii) $\Exp{\mWi} = \mZero$ and $\Exp{\vni} = \vzero$.
        
        (iii) $\Exp{\mPi^T \mWi} = \mZero$ and $\Exp{\mPi^T \vni} = \vzero$.

        (iv) There exist $0 \leq a_0, b_0, c_0 < \infty$ such that
                $\norm{\Exp{\mWi^T \mWi}}_2 \leq a_0$, $\Exp{\vni^T \vni} \leq b_0$, and $\norm{\Exp{\mWi^T \vni}}_2 \leq c_0$, respectively.

\end{assumption}

\begin{remark}\label{rem:comm_ass}
    We do not make any assumptions on the distribution of random quantities involved in \eqref{eq:comm_model}.
    In particular, our model supports systems where the components of $\mPi$ are correlated (such a model may be necessary, for instance, in MIMO systems), or, the components of $\mWi$ and $\vni$ are cross-correlated (a condition necessary for modeling the proposed OTA-C protocol in Section \ref{sec:comm}).
    Note that the models used in previous studies \cite{nedic2009distributed,duchi2011dual,cavalcante_dynamic_2013,michelusi2022non} do not support such systems.
\end{remark}

For proving that our scheme is able to make all agents agree on a common solution to \eqref{eq:relax_set}, we require the following technical assumption on the sequence of subgradients $(\cS(\vhki))_\iN$, and the sequence of outputs $(\gLki)_\iN$ of the local optimization step \eqref{eq:local_scheme_local_opt}:
\begin{assumption} 
    \label{ass:bounded_subg}
    For all $\kN$, the sequences $(\cS(\vhki))_\iN$ and $(\gLki)_\iN$ are bounded a.s.
\end{assumption}
\begin{remark}\label{rem:bounded_subg}
    The assumption of bounded subgradients is standard the literature (e.g.~see \cite{yamada2005adaptive,Slavakis_2006}).
    Additionally, the condition that $(\gLki)_\iN$ is bounded can be satisfied, for example, by using the mapping $\opTki$ as the projection onto a bounded closed convex set. 
    Ensuring the boundedness of $\gLki$ is crucial in practical scenarios, as it allows agents to exchange information over wireless links with finite capacity.
    Note that the conditions stated in Assumption \ref{ass:bounded_subg} are sufficient, but not necessary. 
\end{remark}

The following result establishes, in particular, sufficient conditions for all agents to asymptotically reach consensus on a solution to \eqref{eq:relax_set}.
\begin{theorem}
    \label{thm:main}
    Consider the scheme in \eqref{eq:scheme} in a system where Assumptions \ref{ass:objective}, \ref{ass:local_scheme}, \ref{ass:random_ass} and \ref{ass:bounded_subg} hold.
    Then, the sequence $(\gPi)_\iN$ generated via \eqref{eq:scheme} satisfies the following:

        (i) For any ideal estimate $\bh \in \sUs$, define $\gPs := [(\bh)^T, \dots, (\bh)^T]^T\in\VS$. 
        Then, the sequence $(\norm{\gPi - \gPs}^2)_\iN$ converges, and $(\forall \kN)\ \lim\limits_{i\to\infty} \cT(\vhki) = 0$, a.s.

        (ii) The sequence $((\mI - \mJ) \gPi)_\iN$ satisfies $\liminf\limits_{i\to\infty} \Epp{\norm{(\mI - \mJ) \gPi}} = 0$, a.s.

        (iii) In addition, let the following assumptions hold:
        
            (a) the set $\sUs$ has a nonempty interior w.r.t.~a hyperplane $\Pi\subset\Real^M$, i.e.~$\exists\tilde{\vu}\in\Pi\cap\sUs$ and $\exists\varrho > 0$ such that $\{\vh \in \sM \mid \norm{\vh - \tilde{\vu}} \leq \varrho\} \subset \sUs$, and

            (b) for all $\iN$, the vector $\gZi$ in \eqref{eq:scheme_local_opt} is a function of $\gPi$, i.e.~$(\forall \iN)\ \Epp{\gLi} = \gLi$, a.s.
            
            Then the sequence $(\gPi)_\iN$ converges to a random vector $\hat{\gP} = [\hat{\vh}_1^T, \dots, \hat{\vh}_N^T]^T$ satisfying $\hat{\gP}(\omega) \in \sC$ for a.e.~$\omega\in\Omega$, and $(\forall \kN)\ \lim_{i\to\infty} \cT(\hat{\vh}_k) = 0$, a.s.
        
        (iv) In addition to the assumptions in part (iii) above, let the following hold:
            
            (c) the interior of $\sUs$ is nonempty, i.e.~there exists $\tilde{u}\in\sUs$ and a scalar $\varrho > 0$ such that $\{\vh \in \sM \mid \norm{\vh - \tilde{\vu}} \leq \varrho\} \subset \sUs$, and

            (d) for all $\epsilon>0$, $r>0$, and for a.e.~$\omega \in \Omega$, there exists $\xi(\omega) >0$ such that
            \begin{equation*}
                \begin{aligned}
                    &\inf\limits_{
                        \substack{
                            \Xi(\omega) > \epsilon \\
                            \Pi(\omega) \leq r
                            }
                        } \sum_{k\in\sN} \cT(\vhki(\omega)) \geq \xi(\omega),
                \end{aligned}
            \end{equation*}
            where $\Xi(\omega) := \sum_{k\in\sN} \min\limits_{\vx \in \lev_{\leq 0} \cT} \norm{\vhki(\omega) - \vx}$ and $\Pi(\omega) := \sum_{k\in\sN} \norm{\tilde{\vu} - \vhki(\omega)}$.
            
            Then, for a.e.~$\omega \in \Omega$, the sequence $(\vhki(\omega))_\iN$ converges to the same point $\hat{\vh}(\omega) \in \sUp$ for each agent $k\in\sN$.

    \begin{proof}
        The proof is provided with the supplementary material in Section \ref{apx:thm_main}.
    \end{proof}
\end{theorem}

\begin{remark}\label{rem:random_solution}
    Note that the proposed scheme \eqref{eq:scheme} converges to an $\sUp$-valued random variable.
    Hence, although the optimization problem is deterministic, different runs of \eqref{eq:scheme} may lead to different estimates in the set $\sUp$.
\end{remark}

    
    
    \section{Over-the-Air Computation Protocols}
    \label{sec:comm}

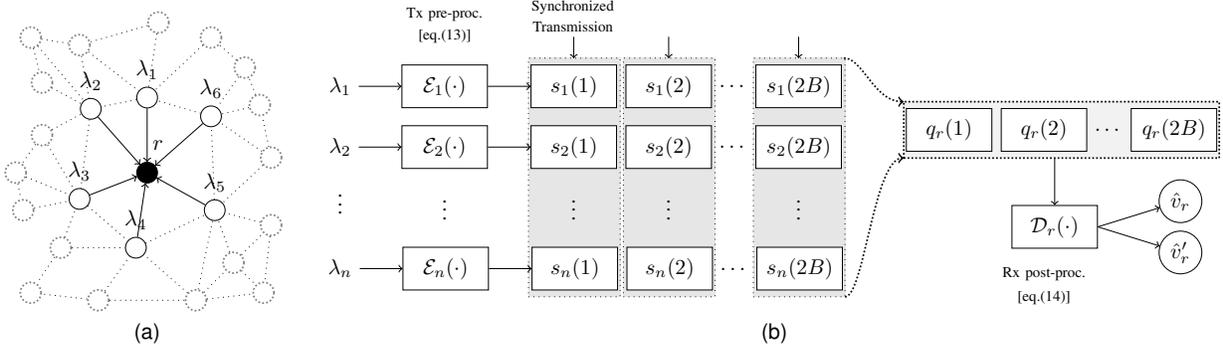
\begin{figure*}[!htb]
    \centering
    \subfloat[]{ 
			\resizebox*{0.23\textwidth}{!}{%
            \begin{tikzpicture}[
                scale=0.6,
                rnode/.style={circle, fill=black, draw=black},
                nnode/.style={circle,fill=none,draw=black},
                wnode/.style={circle,thick,densely dotted,fill=none,draw=black!50},
            ]
            
                \node[rnode, label={[xshift=5pt, yshift=2pt]$r$}] (r) at (0,0) {};
                \node[nnode, label=$\lambda_1$] (l1) at (0.0,2.0) {};
                \node[nnode, label=$\lambda_2$] (l2) at (-1.5,1.7) {};
                \node[nnode, label=$\lambda_3$] (l3) at (-1.8,-0.7) {};
                \node[nnode, label=$\lambda_4$] (l4) at (-0.3,-2.0) {};
                \node[nnode, label=$\lambda_5$] (l5) at (1.8,-1.0) {};
                \node[nnode, label=$\lambda_6$] (l6) at (1.7,1.5) {};
        
                \node[wnode] (w1) at (0.3,4.0) {};
                \node[wnode] (w2) at (-1.8,3.4) {};
                \node[wnode] (w3) at (-2.8,2.4) {};
                \node[wnode] (w4) at (-2.8,1.0) {};
                \node[wnode] (w5) at (-3.5,0.0) {};
                \node[wnode] (w6) at (-3.3,-1.0) {};
                \node[wnode] (w7) at (-2.3,-2.0) {};
                \node[wnode] (w8) at (-1.3,-3.4) {};
                \node[wnode] (w9) at (1.3,-3.4) {};
                \node[wnode] (w10) at (2.3,-2.4) {};
                \node[wnode] (w11) at (3.3,-1.4) {};
                \node[wnode] (w12) at (3.3,0.4) {};
                \node[wnode] (w13) at (3.0,1.8) {};
                \node[wnode] (w14) at (2.6,2.8) {};
                \node[wnode] (w15) at (1.7,3.6) {};
                
                \node[wnode] (w16) at (-3.1,-3.2) {};
                \node[wnode] (w17) at (3.1,-3.2) {};
                \node[wnode] (w18) at (-3.1,3.6) {};
        
                \draw[->] (l1) -- (r);
                \draw[->] (l2) -- (r);
                \draw[->] (l3) -- (r);
                \draw[->] (l4) -- (r);
                \draw[->] (l5) -- (r);
                \draw[->] (l6) -- (r);
        
                \draw[dotted] (w1) -- (l1);
                \draw[dotted] (w15) -- (l1);
                \draw[dotted] (w2) -- (l1);
        
                \draw[dotted] (w2) -- (l2);
                \draw[dotted] (w3) -- (l2);
                \draw[dotted] (w4) -- (l2);
                \draw[dotted] (w18) -- (l2);
                
                \draw[dotted] (w4) -- (l3);
                \draw[dotted] (w5) -- (l3);
                \draw[dotted] (w6) -- (l3);
                \draw[dotted] (w7) -- (l3);
                
                \draw[dotted] (w7) -- (l4);
                \draw[dotted] (w8) -- (l4);
                \draw[dotted] (w9) -- (l4);
                
                \draw[dotted] (w9) -- (l5);
                \draw[dotted] (w10) -- (l5);
                \draw[dotted] (w11) -- (l5);
                \draw[dotted] (w12) -- (l5);
                
                \draw[dotted] (w12) -- (l6);
                \draw[dotted] (w13) -- (l6);
                \draw[dotted] (w14) -- (l6);
                
                \draw[dotted] (w1) -- (w2);
                \draw[dotted] (w2) -- (w18);
                \draw[dotted] (w18) -- (w3);
                \draw[dotted] (w4) -- (w5);
                \draw[dotted] (w5) -- (w6);
                \draw[dotted] (w7) -- (w16);
                \draw[dotted] (w8) -- (w16);
                \draw[dotted] (w8) -- (w9);
                \draw[dotted] (w9) -- (w10);
                \draw[dotted] (w10) -- (w11);
                \draw[dotted] (w9) -- (w17);
                \draw[dotted] (w10) -- (w17);
                \draw[dotted] (w11) -- (w17);
                \draw[dotted] (w12) -- (w13);
                \draw[dotted] (w13) -- (w14);
                \draw[dotted] (w14) -- (w15);
                
                \draw[dotted] (l1) -- (l2);
                \draw[dotted] (l2) -- (l3);
                \draw[dotted] (l3) -- (l4);
                \draw[dotted] (l4) -- (l5);
                \draw[dotted] (l5) -- (l6);
                \draw[dotted] (l6) -- (l1);
                
            \end{tikzpicture}
            }
        }%
        \hfill%
        \subfloat[]{ 
			\resizebox*{0.73\textwidth}{!}{%
            \begin{tikzpicture}[
                blank/.style={draw=none, fill=none},
                box/.style={rectangle, draw=black, fill=white, minimum width=40pt, minimum height=20pt},
                circ/.style={circle, draw=black, fill=none, inner sep=0},
            ]
                \node[blank] (lambda1) at (0,0) {$\lambda_1$};
                \node[box, right=20pt of lambda1.east] (encoder1) {$\mathcal{E}_1(\cdot)$};
                \node[box, right=20pt of encoder1.east] (s11) {$s_1(1)$};
                \node[box, right=4pt of s11.east, anchor=west] (s12) {$s_1(2)$};
                \node[right=0pt of s12.east, anchor=west] (s13) {$\dots$};
                \node[box, right=0pt of s13.east, anchor=west] (s14) {$s_1(2B)$};
                \draw[->] (lambda1.east) -- (encoder1.west);
                \draw[->] (encoder1.east) -- (s11.west);
        
                \node[blank] (lambda2) at (0,-1) {$\lambda_2$};
                \node[box, right=20pt of lambda2.east] (encoder2) {$\mathcal{E}_2(\cdot)$};
                \node[box, right=20pt of encoder2.east] (s21) {$s_2(1)$};
                \node[box, right=4pt of s21.east, anchor=west] (s22) {$s_2(2)$};
                \node[right=0pt of s22.east, anchor=west] (s23) {$\dots$};
                \node[box, right=0pt of s23.east, anchor=west] (s24) {$s_2(2B)$};
                \draw[->] (lambda2.east) -- (encoder2.west);
                \draw[->] (encoder2.east) -- (s21.west);
        
                \node[below=5pt of lambda2.south] (dots1) {$\vdots$};
                \node[below=5pt of encoder2.south] (dots2) {$\vdots$};
                \node[below=5pt of s21.south] (dots3) {$\vdots$};
                \node[below=5pt of s22.south] (dots4) {$\vdots$};
                \node[below=5pt of s24.south] (dots5) {$\vdots$};
        
                \node[blank] (lambdan) at (0,-3) {$\lambda_n$};
                \node[box, right=20pt of lambdan.east] (encodern) {$\mathcal{E}_n(\cdot)$};
                \node[box, right=20pt of encodern.east] (sn1) {$s_n(1)$};
                \node[box, right=4pt of sn1.east, anchor=west] (sn2) {$s_n(2)$};
                \node[right=0pt of sn2.east, anchor=west] (sn3) {$\dots$};
                \node[box, right=0pt of sn3.east, anchor=west] (sn4) {$s_n(2B)$};
                \draw[->] (lambdan.east) -- (encodern.west);
                \draw[->] (encodern.east) -- (sn1.west);
        
                \node[above=17pt of encoder1.north, anchor=south] [xshift=0pt] (Txprep) {{\scriptsize Tx pre-proc.}};
                \node[below=0pt of Txprep.south, anchor=north] [yshift=3pt] (Txprep2) {{\scriptsize  [eq.\eqref{eq:encoder}]}};
        
                \begin{scope}[on background layer]
                    \draw[dotted, draw=black, fill=black!10] ([xshift=-1pt,yshift=3pt]s11.north west) rectangle ([xshift=1pt,yshift=-3pt]sn1.south east);
                    \draw[dotted, draw=black, fill=black!10] ([xshift=-1pt,yshift=3pt]s12.north west) rectangle ([xshift=1pt,yshift=-3pt]sn2.south east);
                    \draw[dotted, draw=black, fill=black!10] ([xshift=-1pt,yshift=3pt]s14.north west) rectangle ([xshift=1pt,yshift=-3pt]sn4.south east);
                \end{scope}
        
                \node[above=20pt of s11.north west, anchor=south west] [xshift=-3pt] (text1) {{\scriptsize Synchronized}};
                \node[below=1pt of text1.south, anchor=north] [xshift=0pt,yshift=5pt] (text12) {{\scriptsize Transmission}};
                \draw[->] ([yshift=13pt]s11.north) -- ([yshift=3pt]s11.north);
                \draw[->] ([yshift=13pt]s12.north) -- ([yshift=3pt]s12.north);
                \draw[->] ([yshift=13pt]s14.north) -- ([yshift=3pt]s14.north);

                \node[box, right=30pt of s14.south east, anchor=north west] (q1) {$q_r(1)$};
                \node[box, right=4pt of q1.east, anchor=west] (q2) {$q_r(2)$};
                \node[right=0pt of q2.east, anchor=west] (q3) {$\dots$};
                \node[box, right=0pt of q3.east, anchor=west] (q4) {$q_r(2B)$};
        
                \begin{scope}[on background layer]
                    \draw[thick, densely dotted, draw=black, fill=black!5] ([xshift=-1pt,yshift=3pt]q1.north west) rectangle ([xshift=1pt,yshift=-3pt]q4.south east);
                \end{scope}
        
                \draw [densely dotted, thick, ->] ([xshift=1pt, yshift=3pt]s14.north east) .. controls ([xshift=10pt, yshift=3pt]s14.north east) and ([xshift=-15pt, yshift=3pt]q1.north west) .. ([xshift=-1pt, yshift=3pt]q1.north west);
                \draw [densely dotted, thick, ->] ([xshift=1pt, yshift=-3pt]sn4.south east) .. controls ([xshift=10pt, yshift=-3pt]sn4.south east) and ([xshift=-15pt, yshift=-3pt]q1.south west) .. ([xshift=-1pt, yshift=-3pt]q1.south west);

                \node[box, below=25pt of q2.south, anchor=north] [xshift=5pt] (decoder) {$\mathcal{D}_r(\cdot)$};

                \draw[->] ([xshift=5pt,yshift=-3pt]q2.south) -- (decoder.north);
        
                \node[circ, minimum size=20pt, right=of decoder.east] [yshift=12pt] (v1) {$\est{v}_r$};
                \node[circ, minimum size=20pt, right=of decoder.east] [yshift=-12pt] (v2) {$\est{v}'_r$};
        
                \node[below=5pt of decoder.south, anchor=north] [xshift=-5pt] (Rxproc) {{\scriptsize Rx post-proc.}};
                \node[below=-2pt of Rxproc.south, anchor=north] (Rxproc2) {{\scriptsize [eq.\eqref{eq:r_postproc}]}};
                
                \draw[->] (decoder.east) -- (v1);
                \draw[->] (decoder.east) -- (v2);
                
            \end{tikzpicture}
            }
        }
    \caption{
        The proposed OTA-C protocol. 
        Figure (a) shows the subgraph of the agent $r$ and its neighbors in $\sN_r$.
        Figure (b) gives an overview of the protocol. 
        Every agent $j\in\sN_r$ encodes its information $\lambda_j$ into $2B$ symbols using \eqref{eq:encoder}.
        Then, all agents in $\sN_r$ transmit synchronously symbols $(s_j(b))_{j\in\sN_r}$, for every $\bB$, over a noisy and fading WMAC.
        In response, the agent $r$ receives $2B$ symbols, which are then post-processed using \eqref{eq:r_postproc} to obtain estimates of $v_r$ and $v'_r$.
    }
    \label{fig:OTAC-schematics}
\end{figure*}

In this section, we present a novel communication protocol to implement the consensus step \eqref{eq:scheme_consensus} using the OTA-C technology \cite{goldenbaum2013robust,frey_ota_2019}.
To motivate our design, as an example, consider a multi-agent system represented by an undirected and time-invariant graph $\sG := (\sN, \sE)$ with weights $(\forall (j,r)\in\sE)\ \rz \in\Realp$.
Suppose that each agent $\kN$ possesses a scalar $\lambda_{k,0}$ at time $i=0$, and their objective is to compute the average $\bar{\lambda}_0 := (1/N)\sum_{\kN} \lambda_{k,0}$, $N:=\abs{\sN}$, via local exchanges over the network graph $\sG$.
Towards this goal, each agent $r\in\sN$ runs the following first-order consensus protocol that generates a sequence $(\lambda_{r,i})_\iN$ via:
\begin{align}\label{eq:lin_con}
     \lambda_{r,i+1} = 
        \Big(1 - \beta_i \sum_{j\in\sN_r} \rz\Big) \lambda_{r,i} + \beta_i \sum_{j\in\sN_r} \rz\ \lambda_{j,i},
\end{align}
where the set of neighbors of agent $r$ is defined as $\sN_r := \{j: (j,r)\in\sE\}$, and $(\beta_i)_\iN$ is a sequence of step-sizes in $(0, 1]$.
For each agent $r\in\sN$, under certain conditions (see \cite{kar2008distributed} for details), the sequence $(\lambda_{r,i})_\iN$ converges asymptotically to the average $\bar{\lambda}_0$ \cite{kar2008distributed,fagnani2008randomized,pescosolido2008average}.
Note that, to implement the protocol \eqref{eq:lin_con}, agent $r$ only needs the scalars $v_{r,i} := \sum_{j\in\sN_r} \rz \lambda_{j,i}$, and $v'_{r,i} := \sum_{j\in\sN_r} \rz$, while the explicit weights $(j\in\sN_r)\ \rz$ need not be known.
In the consensus step of our algorithm, agents update their local estimate $\vhki\in\Real^M$ using multiple instances (one for each coordinate of $\vhki$) of a protocol similar to \eqref{eq:lin_con}, given later in \eqref{eq:agents_protocol}.
The proposed OTA-C protocol is designed to enable every agent $r\in\sN$ to obtain estimates of the scalars $v_{r,i}$ and $v'_{r,i}$ in a truly decentralized manner over a time-varying network with fading and noisy wireless channels.

With this aim, in Section \ref{sec:OTAC_design}, we first describe the OTA-C protocol design that implements consensus step for scalar-valued inputs and outputs.
The protocol is extended in Section \ref{sec:deploy_OTAC} for exchanging vector-valued information to implement \eqref{eq:local_scheme_consensus}.
The communication mechanism of the resulting protocol, given in \eqref{eq:agents_protocol}, follows the abstract model $\opLi$ in \eqref{eq:comm_model}.
In the proposed algorithm for optimization, each agent iterates between the local optimization step \eqref{eq:local_scheme_local_opt} and the proposed OTA-C based consensus step \eqref{eq:agents_protocol}.
In Proposition \ref{prop:OTA_scheme}, we prove that the sufficient conditions for convergence in Theorem \ref{thm:main} hold for the proposed OTA-C based algorithm.

\subsection{The OTA-C protocol design} \label{sec:OTAC_design}
In the following, the time index $i$ is omitted, as the protocol remains the same for every $\iN$.
Consider the subgraph of a (random) graph $\sG:=(\sN, \sE)$ formed by agent $r\in\sN$ and its neighbors in $\sN_r := \{j\in\sN \mid (j, r) \in \sE\}$, as shown in Figure \ref{fig:OTAC-schematics}(a).
The simultaneous transmissions of wireless signals by agents in $\sN_r$, and reception by agent $r$, is described by the wireless multiple access channel (WMAC) model~\cite{goldenbaum2013robust,frey_ota_2019}.
As described later in Section \ref{sec:preproc}, a single run of the OTA-C protocol requires $2B$ uses of the WMAC, where, in any $b$th use, for $b=1,\dots, 2B$, the WMAC model takes as input the complex-valued symbols $(\forall j\in\sN_r)\ s_j(b)$ from the neighbors of agent $r$,
and outputs a random number $q_r(b)\in\Complex$ given by: $(\forall r\in\sN)(b=1,\dots,2B)$
\begin{align}
    \label{eq:WMAC}
    q_r(b) := \sum_{j\in\sN_r} \xi_{jr}(b)\, s_j + w_r(b),
\end{align}
where 
$\xi_{jr}$ and $w_r$ are complex-valued random variables representing the channel fading between agents $j$ and $r$, and the receiver noise at agent $r$, respectively.
We make the following assumptions on the WMAC model.
\begin{assumption}\label{ass:WMAC}
    For all agents $(r,l)\in\sN\times\sN$, and their neighbors $(\forall j\in\sN_r)$ and $(\forall k\in\sN_l)$, the following statements hold: \hspace{\bulletgap}
    (i) $\Exp{w_r} = 0$. \hspace{\bulletgap} 
    (ii) $\Exp{\abs{w_r}^2} < \infty$. \hspace{\bulletgap} 
    (iii) $\Exp{\abs{w_r}^2 \abs{w_l}^2} < \infty$. \hspace{\bulletgap}
    (iv) $\Exp{\abs{\xi_{jr}}^2} < \infty$. \hspace{\bulletgap} 
    (v) $\Exp{\abs{\xi_{jr}}^2 \abs{\xi_{lk}}^2} < \infty$. \hspace{\bulletgap}
    (vi) The random variables $\xi_{jr}$ and $w_l$ are independent. \hspace{\bulletgap}
    (vii) $\xi_{jr}(b) = \xi_{jr}(b+1)$ for every odd $b$, where $b\in\Natural$, $B>0$.
\end{assumption}

\begin{remark} \label{rem:sync}
    The OTA-C only requires a coarse (frame-level) synchronization of agents in practice \cite{goldenbaum2013robust}. 
    This property enables all agents in a decentralized network to transmit and receive simultaneously.
\end{remark}

\begin{remark}\label{rem:WMAC_ass}
    The previous studies \cite{goldenbaum2013robust,frey_ota_2019} impose stronger conditions on the WMAC compared to Assumptions \ref{ass:WMAC}(i)-(vi). For example, \cite{goldenbaum2013robust} assumes i.i.d.~channels, and in \cite{frey_ota_2019}, all channels are assumed to have a unit variance.
    Furthermore, Assumption \ref{ass:WMAC}(vii) ensures that both realizations of the random numbers $(\forall j\in\sN_r)\ \rz$ corresponding to $v_r$ and $v'_r$ take the same values (see \eqref{eq:lin_con} and the subsequent discussion). 
    For instance, such assumption is valid in an OFDM system with coherence bandwidth spanning at least two sub-carriers.
\end{remark}

To obtain estimates of $v_r$ and $v'_r$, every transmitting agent in $\sN_r$ and the receiving agent $r$ perform pre- and post-processing operations, respectively, as described in Sections \ref{sec:preproc} and \ref{sec:postproc}.
Note that, unlike previous studies \cite{goldenbaum2013robust,abari2016over} that require information about individual weights $(j\in\sN_r)\ \rnu_{jr}$, our protocol only requires an estimate of their sum, which it obtains using the same OTA-C protocol by transmitting a second (dummy) scalar value, as described in the following.

\subsubsection{Transmit pre-processing} \label{sec:preproc}
Let the scalar-valued input of any agent $j\in\sN_r$ be $\lambda_j$, where $(\forall \kN)$ $\lambda_k \in \sS:= [\dmin, \dmax] \subset \Real$, $\dmin < \dmax$, and $(\dmin,\dmax)$ are assumed to be known to all agents.
Every transmitter $j\in\sN_r$ encodes $\lambda_j$ into $2B$ complex symbols $(s_j(1), \dots, s_j(2B))$, given by: 
\begin{equation}
    \begin{aligned}
        \label{eq:encoder}
        s_j(b) := \begin{cases}
            \sqrt{g_j(\lambda_j)}\ U_j(b) & \text{if $b$ is odd} \\
            \sqrt{g_j(\dmax)}\ U_j(b) & \text{if $b$ is even},
        \end{cases}
    \end{aligned}
\end{equation}
where $g_j(x) := (P_j/\Delta) (x - \dmin)$, $\Delta := \dmax - \dmin$, ensures the peak transmit power constraint $P_j>0$ at agent $j$,\footnote{Later, in Assumption \ref{ass:network}(i), we use the same peak transmit power for all agents to realize an undirected graph (also see Remark \ref{rem:network}).}~%
and the complex scalars $(U_j(1),\dots,U_j(2B))$ are i.i.d.~random variables such that $(\forall \bB)$ $\Exp{U_j(b)} = 0$ and $\abs{U_j(b)} = 1$.\footnotemark~%
\footnotetext{The scalar $U_j(b)\in\Complex$ can be sampled, e.g., uniformly from $\{-1, 1\}$ \cite{frey_ota_2019}, or uniformly on unit circle in complex plane \cite{goldenbaum2013robust}.}
For each $b=1,\dots,2B$, every agent $j\in\sN_r$ transmits symbol $s_j(b)$ synchronously to agent $r$ over the WMAC \eqref{eq:WMAC} (see Remark \ref{rem:sync}).
Thus, every run of the OTA-C protocol requires $2B$ WMAC resources.

\subsubsection{Receiver post-processing}\label{sec:postproc}
In response, agent $r$ receives $2B$ symbols $(q_r(1), \dots, q_r(2B))$, where $q_r(b)$, for $b=1,\dots,2B$, is given by \eqref{eq:WMAC}.
Let $\vy_r := [q_r(1), q_r(3), \dots, q_r(2B-1)]^T$ and $\vy'_r := [q_r(2), q_r(4), \dots, q_r(2B)]^T$.
Agent $r$ obtains the estimates of $v_r$ and $v'_r$ by evaluating:
\begin{equation}
\label{eq:r_postproc}
    \begin{aligned}
        \est{v}_r &:= \Delta (B^{-1}\norm{\vy_r}^2 - \eNoise) + \dmin \est{v}'_r \\ \est{v}'_{r} &:= B^{-1}\norm{\vy'_r}^2 - \eNoise,
    \end{aligned}
\end{equation}
where the receiver noise power $\eNoise$, defined as $\eNoise := \Exp{\abs{w_r}^2}$, is assumed to be known to agent $r$.
In Lemma \ref{lem:unbiased_est} below, we establish that the scalars $\est{v}_{r}$ and $\est{v}'_{r}$ in \eqref{eq:r_postproc} are unbiased estimates of $v_r$ and $v'_r$, respectively, conditioned on the random variables $(\forall j\in\sN_r)\ \rz$.

\begin{lemma}\label{lem:unbiased_est}
    Consider the subgraph of $\sG$ formed by any agent $r\in\sN$ and its neighbors in $\sN_r$, as shown in Figure \ref{fig:OTAC-schematics}(a).
    Let Assumption \ref{ass:WMAC} hold for the WMAC model \eqref{eq:WMAC}, and suppose that agent $r$ and its neighbors in $\sN_r$ implement the post-~and pre-processing operations described above, respectively.
    Then, for all $r\in\sN$, the scalars $\est{v}_{r}$ and $\est{v}'_{r}$ obtained by agent $r$ in \eqref{eq:r_postproc} can be expressed as:
    \begin{align}\label{eq:otac_scalar}
        \est{v}_{r} = \sum_{j\in\sN_r} \rnu_{jr} \lambda_j + \eta_r, \quad 
        \est{v}'_{r} = \sum_{j\in\sN_r} \rnu_{jr} + \eta'_r,
    \end{align}
    where $(\forall j\in\sN_r)\ \rnu_{jr} := (P_j/B) \sum_{b=1}^B \abs{\xi_{jr}(2b-1)}^2$, and random variables $\eta_r$ and $\eta'_r$ are zero-mean.
    More precisely, the conditional expectation of the absolute error in estimation is zero, i.e.~$\Exp{\abs{\est{v}_{r} - v_r} \cond (\rz)_{j\in\sN_r}} = 0$ and $\Exp{\abs{\est{v}'_{r} - v'_r} \cond (\rz)_{j\in\sN_r}} = 0$.
    \begin{proof}
        The proof is provided with the supplementary material in Section \ref{apx:lem_otac}.
    \end{proof}
\end{lemma}


\subsection{Deploying the proposed scheme using OTA-C}\label{sec:deploy_OTAC}
In this section, the OTA-C protocol proposed in Section \ref{sec:OTAC_design} is implemented separately for each scalar coefficient of the $M$ dimensional vector inputs $(\gL_{j,i})_{j\in\sN_r,\iN}$ in the consensus step \eqref{eq:local_scheme_consensus}.
We reinstate the indexing $\iN$ in the following to identify the iterates.
As a result, the information obtained by agent $r$ can be expressed as: $(\forall\iN)(\forall r\in\sN)$
\begin{equation}
    \label{eq:OTA-C-output}
    \begin{aligned}
        \est{\vv}_{r,i} &= \sum_{j\in\sN_r} \mZ_{jr,i}\ \gL_{j,i} + \ovn_{r,i}, \\
        \est{\vv}'_{r,i} &= \sum_{j\in\sN_r} \mZ_{jr,i}\ \vOne_M + \vnn_{r,i},
    \end{aligned}
\end{equation}
where $\mZ_{jr,i}:=\diagm([\rnu_{jr,i}^{(1)}, \dots, \rnu_{jr,i}^{(M)}]^T)$ is a diagonal matrix, and $\vOne_M\in\Real^M$ is the vector of all ones.
The scalars $\rnu_{jr,i}^{(m)}\in\Realp$, for all $m=1,\dots,M$, (defined explicitly in Lemma \ref{lem:unbiased_est} above) are functions of realizations of the channel for the $m$th run of the OTA-C at time $i$.
Similarly, vectors $\ovn_{r,i}$ and $\vnn_{r,i}$ are formed by stacking $M$ samples of random variables $\eta_{r,i}$ and $\eta'_{r,i}$, respectively.

Next, we describe the iterative scheme implemented by agent $r\in\sN$, using $\est{\vv}_{r,i}$ and $\est{\vv}'_{r,i}$ which have been obtained in \eqref{eq:OTA-C-output} via the OTA-C.
Starting with an arbitrary $\vh_{r, 0} \in \Real^M$, agent $r$ implements the local optimization step \eqref{eq:local_scheme_local_opt}, yielding a vector $\gL_{r,i}$ at time $i$.
In the consensus step, agent $r$ implements: $(\forall \iN)(\forall r\in\sN)$
\begin{align}
    \label{eq:agents_protocol}
    \vh_{r, i+1} = (\vOne_M - \beta_i\gamma_i\ \est{\vv}'_{r,i}) \odot \gL_{r,i} + \beta_i\gamma_i\ \est{\vv}_{r,i},
\end{align}
where $\odot$ stands for element-wise product.
The scalar $\gamma_i$ is designed such that each component of $\Exp{(\vOne_M - \beta_i\gamma_i\vv'_{r,i})}$ is nonnegative for every $r\in\sN$.
More precisely: $(\forall \iN)$
\begin{align}
    \label{eq:gamma}
    0 < \gamma_i \leq \big(\max_{r\in\sN} \sum_{j\in\sN_r} \Exp{\rni^{(m)}}\big)^{-1} < \infty,
\end{align}
where $\Exp{\rni^{(m)}} = (P_k/B) \sum_{b=1}^B \Exp{\abs{\xi_{jr,i}^{(m)}(2b-1)}^2}$. 
Note that each agent in the network must employ the same sequences $(\gamma_i)_\iN$ and $(\beta_i)_\iN$, and the product $\gamma_i\beta_i$ leads to a trade-off between information diffusion and noise suppression.

By defining $\ogP_i$ and $\ogL_i$ as vectors formed by stacking $\vh_{r,i}$ and $\gL_{r, i}$ for all $r\in\sN$, respectively, 
\eqref{eq:agents_protocol} becomes
\begin{align}\label{eq:ota_protocol}
    (\forall  \iN) \qquad \ogP_{i+1} = (1 - \beta_i)\,\ogL_i + \beta_i\, \check{\opL}_i(\ogL_i),
\end{align}
where $\check{\opL}_i(\gP) := (\oPi + \oWi)\, \gP\, +\, \ovn_i$ describes the information exchange in the network under the proposed OTA-C protocol.
The matrix $\oPi$ is a block-matrix with $[\oPi]_{(p,q)}\in\Real^{M\times M}$ for each $(p,q)\in\sN\times\sN$ given by:
\begin{align} \label{eq:matrix_P}
    [\oPi]_{(p,q)} := \begin{cases}
        \mI_M - \gamma_i \sum_{j\in\sN_p} \mZ_{jp, i} &\text{ if } p = q, \\
        \gamma_i \mZ_{qp, i} &\text{ if } (p,q) \in \sE_i, \\
        \mZero &\text{ otherwise},
    \end{cases}
\end{align}
where $\mI_M\in\Real^{M\times M}$ is the identity matrix of dimension $M$.
Moreover, $\oWi := \diagm\left([\vnn_{1, i}^T, \dots, \vnn_{N, i}^T]^T\right)$, and $\ovni := [\ovn_{1, i}^T, \dots, \ovn_{N, i}^T]^T$,
where $\vnn_{k,i}$ and $\ovn_{k,i}$ are defined in \eqref{eq:OTA-C-output}.

\begin{assumption}\label{ass:network}
    At any time $\iN$, given the random graph $\sGi := (\sN, \sE_i)$, the following statements hold: \hspace{\bulletgap}
    (i) The expected graph $\sGb$ (defined in Section \ref{sec:prelims}) is \textit{connected}.  \hspace{\bulletgap}
    (ii) For all $(j,r)\in\sEi$, the edge weights of $\sGb$ are reciprocal, i.e.~$(\forall (j,r)\in\sE)\ \Exp{\rni} = \Exp{\rnu_{rj,i}}$, and hence, the graph $\sGb$ is \textit{undirected}.
    \end{assumption}

\begin{remark}\label{rem:network}
    Assumption \ref{ass:network}(i) and (ii) are standard in the literature dealing in graphs with time-varying topologies \cite{olfati2007consensus,nedic2009distributed,cavalcante_dynamic_2013}.
    Condition \ref{ass:network}(ii) ensures that the matrix $\mPi$ in \eqref{eq:comm_model} resulting from the random graph $\sGi$ satisfies the conditions in Definition \ref{def:consensus_matrix}.
    In practice, such a graph can be realized simply by adopting equal peak transmit powers for every agent i.e.~$P_p = P_q$ for every $(p,q) \in \sE$.
\end{remark}
        

\begin{remark}\label{rem:OTAC}
    At any iteration $\iN$, the total number of WMAC resources (channel uses) required to run the proposed OTA-C protocol over the entire network is equal to $2BM$.
    Hence, the parameter $B$ provides a trade-off between accuracy and resource usage.
    In contrast, for channel separation-based communication schemes, such as gossip or broadcast, the required amount of wireless resources grow with the number of agents in the network.
    In Section \ref{sec:sim}, we show that the proposed OTA-C protocol outperforms the separation-based schemes using an equivalent amount of wireless resources.
\end{remark}

The proposed scheme iterates between the local optimization step \eqref{eq:local_scheme_local_opt}, and the OTA-C based implementation of the consensus step \eqref{eq:ota_protocol}.
To prove that the scheme converges asymptotically to an optimal solution, we need to prove that the set of assumptions in Theorem \ref{thm:main} are satisfied.
In this direction, note that the conditions in Assumptions \ref{ass:objective} and \ref{ass:bounded_subg} are already satisfied for the given optimization problem \eqref{eq:relax_set}, and Assumption \ref{ass:local_scheme} holds for the system by design.
In Proposition \ref{prop:OTA_scheme} below, we prove that the remaining conditions in Theorem \ref{thm:main}, i.e.~conditions in Assumption \ref{ass:random_ass}, are satisfied for the proposed OTA-C based scheme, provided that Assumption \ref{ass:WMAC} holds.
\begin{proposition} \label{prop:OTA_scheme}
    Consider a wireless network where agents exchange information using the proposed OTA-C protocol, and let Assumption \ref{ass:WMAC} hold.
    Then, for the model $\check{\opL}_i$ in \eqref{eq:ota_protocol}, the conditions in Assumption \ref{ass:random_ass} are satisfied.
    \begin{proof}
        The proof is provided with the supplementary material in Section \ref{apx:prop_ota}.
    \end{proof}
\end{proposition}

    \section{Exemplary Application}
    \label{sec:sim}



In this section, we demonstrate the efficacy of the proposed algorithm \eqref{eq:scheme} by applying it to an online supervised learning problem in a distributed setting.
Specifically, the data is made available sporadically to a decentralized and time-varying network of sensors.
The dataset used for simulation consists of real-world ocean temperature measurements taken at various locations and depths in the Gulf of Mexico \cite{nceiDataset}.
The system consists of $N$ agents (sensors), indexed by $\sN:=\{1, \dots, N\}$, and they exchange information over a wireless network with fading and noisy channels. 
All agents periodically move to new locations to obtain new temperature measurements, which result in a change of the network graph.
In Section \ref{sec:sim_optprob}, we formulate the optimization problem, and present the distributed algorithm based on \eqref{eq:scheme}. 
The problem is modeled using the dictionary-based multikernel approach\cite{yukawa2012multikernel} with Random Fourier Features (RFF) approximation of the kernels \cite{shen2020multikernel}.
Furthermore, a sequence of bounded perturbations $(\zeta_i \vzki)_\iN$ is designed for the local optimization step \eqref{eq:scheme_local_opt} to promote sparsity in $\gLki$. 
The sparsity in vector $\gLki$ is beneficial as it leads to significant energy savings during the implementation of the consensus step \eqref{eq:local_scheme_consensus} using the proposed OTA-C protocol.
In Section \ref{sec:sim_net}, we discuss implementation details related to the realization of the wireless network, and the parameters used for simulation.
The results are presented in Section \ref{sec:sim_results}, where we analyze the proposed distributed scheme for various parameters and settings, and compare its performance with other standard methods.


\subsection{Optimization Problem and Algorithm}\label{sec:sim_optprob}
Let function $f:\sD\to [0,1]$ maps any location $\vx\in \sD := [0,1000]^3$ in the ocean to its temperature value $f(\vx)\in [0, 1]$.\footnote{The locations and temperatures in the original dataset are linearly scaled to $[0, 1000]^3$ and $[0,1]$, respectively.}~%
The distances are measured in meters.
In the following, we describe the approach used for modeling the function $f$.
Let $\{\Ker_1, \dots, \Ker_{2L}\}$ be a set of $2L$ distinct Gaussian kernels,
where for $l = 1,\dots,L$, $\Ker_l:\sD\times\sD\to\Realp$, and for $l = L+1,\dots, 2L$, $\Ker_l= -\Ker_{l-L}$.%
\footnote{Each kernel is included in the dictionary twice, once with a positive sign and once with a negative sign. As a result, the resulting model consists of nonnegative parameters $\vhki$ (model weights). This property is particularly advantageous under the proposed OTA-C protocol, as it allows for the omission of any components in the vector $\gLki$ that have a value of zero during transmission, and hence, saving energy.}~%
Moreover, each kernel is represented (approximately) using $P$ RFFs \cite{Rahimi2007,shen2020multikernel}.
Under this framework, the function $f$ is modeled as\cite{shin2018distributed,shen2020multikernel}:
\begin{align} \label{eq:est_model}
    \hf(\vx) := 
                \sum_{l=1}^{2L} \sum_{p=1}^P h_{l, p}\  \pv_{l,p}(\vx) = \vh^T\ \pvv(\vx),
\end{align}
where, for all $l=1,\dots,2L$ and $p=1,\dots,P$, the nonnegative scalars $(h_{l, p})$ are the parameter to be learned.
For each $l=1,\dots,L$, define $\pv_l$ : $\Omega \times \sD \to \Real$ : $(\omega,\vx) \mapsto \sqrt{2/P}\cos(\vx^T\ow_l(\omega) + \pb_l(\omega))$, where random variables $\ow_l$ and $\pb_l$ are Gaussian and uniformly distributed, respectively, such that the relation $\mathbb{E}[\pvv_l(\vx)\pvv_l(\vy)] = \Ker_l(\vx, \vy)$ holds \cite{Rahimi2007}.
The scalars $(\pv_{l,p}(\vx))_{p=1}^P$ are i.i.d.~samples of $\pv_l(\vx)$.
Note that $M=2LP$ for $\vh\in\Real^M$ in \eqref{eq:est_model}.

The cost function $\cT:\Real^{M} \to \Realp$ is given by \cite{yukawa2012multikernel}: $(\forall \iN)(\forall \kN)\ \cT(\vh) := \norm{\vh - \opP_{k,i}(\vh)} \norm{\vhki - \opP_{k,i}(\vhki)}$, where $\vhki$ is the estimate of agent $\kN$ at time $\iN$, and operator $\opP_{k,i}$ is the projection operator onto the set $\sQki := \{\vh : \abs{\vh^T \pvv(\vx_{k,i}) - \ey_{k,i}} \leq \epsilon_k\}$ (see \cite{yukawa2012multikernel} for the expression of $\opP_{k,i}$).
The scalar $\epsilon_k$ is a design parameter, chosen such that Assumption \ref{ass:objective}(ii) is satisfied with high probability.
Furthermore, for all $\iN$ and $\kN$, we define the operator $\opTki$ as the orthogonal projection onto the set $\sXki$.
As a result, each component of the input $\gLki$ is restricted to $\sS := [\dmin, \dmax] = [0, 1]$.

By defining $\vy_{k,i} := \vhki - \rmuki (\vhki - \opP_{k,i}(\vhki))$, the local optimization step \eqref{eq:local_scheme_local_opt} can be written as $\gLki = \opTki(\vy_{k,i} + \zeta_i \vzki)$.
To promote sparsity in $\gLki$, we design the sequence of bounded perturbations $(\zeta_i \vz_{k,i})_\iN$ such that the resulting vector iteratively reduces the reweighted $\ell_1$-norm of vector $\vy_{k,i}$.\footnote{See \cite{candes2008enhancing} for details on reweighted $\ell_1$-norm for enhansing sparsity.}~%
More precisely, the vector $\vzki$ is designed as: $(\forall \iN)(\forall \kN)$
\begin{align} \label{eq:perturb_seq}
     \vzki := \zeta_i^{-1}\ \left(\fq_{k,i}(\vy_{k,i}) - \vy_{k,i} \right),
\end{align}
where $(\forall\iN)\ \zeta_i>0$ and the series $\sum_\iN \zeta_i$ is convergent.
The operator $\fq_{k,i}:\Real^M\to\Real^M$ is defined element-wise for each $m=1,\dots,M$ as
$\fq_{k,i}^{(m)}(x) := \text{sign}(x)\ \left[\abs{x} - \zeta_i(\abs{\vy_{k,i-1}[m]} + \varsigma)^{-1}\right]_+$,
where $\text{sign}(a) := a/\abs{a}$, and $[a]_+ := \max(a, 0)$.

The following proposition verifies that $(\zeta_i\vzki)_\iN$ is a sequence of bounded perturbations.
\begin{proposition} \label{prop:perturb}
    Consider the sequences $(\zeta_i)_\iN$ and $(\vzki)_\iN$ as defined in \eqref{eq:perturb_seq}, for all $\kN$.
    Then, for every $\kN$, the sequence $(\zeta_i\vzki)_\iN$ is a \emph{sequence of bounded perturbations} in the sense of Definition \ref{def:perturb}.
    \begin{proof}
        The proof is provided with the supplementary material in Section \ref{apx:prop_sim}.
    \end{proof}
\end{proposition}

In summary, the scheme implemented by any agent $\kN$ at time $\iN$ is given by:
\begin{align*}
    \gLki           &= \opTki[\fq_{k,i}(\vhki - \rmuki (\vhki - \opP_{k,i}(\vhki)))]   \\
    \vh_{k, i+1}    &=  (\vOne - \beta_i\gamma_i\ \vv'_{k,i}) \odot \gLki + \beta_i\gamma_i\ \vv_{k,i},
\end{align*}
where $\vv_{k,i}$ and $\vv'_{k,i}$ are obtained via the proposed OTA-C protocol in \eqref{eq:OTA-C-output}, and the value of the scalar $\gamma_i$ is chosen based on \eqref{eq:gamma}.


\subsection{Implementation details}\label{sec:sim_net}
We simulate a network of agents in $\sN$, given by graph $\sGi:=(\sN, \sE_i)$, where $\sE_i = \sN \times \sN$ for all $\iN$, i.e.~the graph is complete.
At random intervals of time $\iN$ (Poisson distributed with mean duration $300$), agents in $\sN$ move to a new position sampled randomly from the training dataset, and obtain a noisy temperature measurement.
Let $\ssT\subset\Natural$ denotes the set of time indices at which a new measurement is obtained. 
Then, the measurement obtained by any agent $\kN$ at time $\tau\in\ssT$ is $\ey_{k,\tau} := f(\vx_{k,\tau}) + n_{k,\tau}$, where measurement noise $(\forall\kN)(\tau\in\ssT)$ $n_{k,\tau}\in\Real$ is Gaussian with zero mean and variance $0.09$.
Agents are equipped with half-duplex transceivers, and at every $\iN$, each of them flips a coin to decide whether to transmit or receive in the time slot $\iN$.
In the WMAC model \eqref{eq:WMAC}, for all $\iN$ and $r\in\sN$, we use zero-mean circularly symmetric Complex Gaussian distributions to model noise $w_{r,i}$ and fading channels $(\forall j\in\sN_r)\ \xi_{jr,i}$ with variances $\sigma_w^2$ and $\sigma_{jr,i}^2$, respectively.
Furthermore, we assume (in addition to Assumption \ref{ass:WMAC}) that for every pair $(j,r) \in \sEi$, and for every $b=1, \dots, 2B$, the channels corresponding to each run of OTA-C for $m=1, \dots, M$ are coherent, i.e.~$\xi_{jr,i}^{(m)}(b) = \xi_{jr,i}^{(m')}(b)$ for all $m,m'\in\{1, \dots, M\}$.
Recall from Section \ref{sec:OTAC_design} that the transmit symbols $(s_j(b))_{j\in\sN_r}$ corresponding to the even $b$ in \eqref{eq:encoder} are used for obtaining an estimate of $v'_r=\sum_{j\in\sN_r} \rnu_{jr}$, where $\rnu_{jr}$ is a function of realizations of the channel $\xi_{jr}$.
Since the channels for each $m=1,\dots,M$ are coherent, we only transmit symbols $(s_j(b))_{j\in\sN_r}$ corresponding to even $b$ once for $m=1$, resulting in a total of $(M+1)B$ symbols transmitted in every OTA-C run.

The path loss follows the Friis' formula, i.e.~channel power between any two agents $j$ and $r$ at time $i$ is given by $\sigma_{jr,i}^2 = \lambda^2/(4\pi d_{jr,i})^2 P_T$, where $d_{jr,i}$ is distance between them, $\lambda$ is the wavelength for $3$GHz carrier, and $P_T=0$dBm is the transmit power, assumed to be the same for all the transmitters.
The receiver noise power is the same for all devices, and it is specified as the signal-to-noise-ratio (SNR) observed at a distance $500$m from a transmitter. 
The bandwidth of the signal is $1$MHz.
The step-size $(\forall \iN)\ \beta_i := (\lfloor i/50 \rfloor)^{-0.51}$.
In \eqref{eq:perturb_seq}, the scalars $(\forall \iN)$ $\zeta_i$ and $\varsigma_i$ are designed as $\zeta_i := 10^{-7} (\lfloor i / 100 \rfloor + 1)^{-1}$, and $\varsigma_i := \max(25 \zeta_i, 10^{-11})$.
Other parameters are $L=2$, $P=25$, $B=10$, $\epsilon_k = 0.3$ and $\rmuki = 0.5$.

We compare the performance of the proposed scheme with alternative approaches that differ solely in the communication protocols employed to execute the consensus step:
(i) \emph{NOC}: no communication scheme, where agents only implement the local optimization step \eqref{eq:local_scheme_local_opt}, 
(ii) \emph{CEN}: the ideal centralized scheme, where each agent obtains the exact (noiseless) average of the parameters from all agents in the network at every step,
(iii) \emph{DBC}: a digital broadcast scheme based on the channel separation principle with prior scheduling of agents, and
(iv) \emph{ANB}: an analog broadcast scheme, similar to DBC, but exchanging analog values modulated on signal power.
Note that, in schemes (i)--(iv) above, and for the vanilla version of the proposed scheme, denoted by \emph{OTA-C}, the local optimization step is implemented without sparsity promoting perturbations, i.e.~$\zeta_i=0$ for all $\iN$. In this case, the parameter $\vhki$ can take values from the set $\sS := (\dmin,\dmax)=(-1,1)$, since the weights can take negative values, and hence, instead of using $2L$ kernels in \eqref{eq:est_model}, we use only $L$ kernels.
The variant of the proposed OTA-C scheme with sparsity promoting iterates, as described in Section \ref{sec:sim_optprob}, is denoted by \emph{OTA-CS}.

The NOC and CEN schemes are benchmark schemes for the worst and the best communication protocols, respectively.
In the DBC scheme, assuming Rayleigh fading channels, the rate of digital communication is selected such that it amounts to less than $20\%$ outage probability at a distance of $500$m for a given SNR.
The $M$ dimensional input vector $\gLki$ is encoded into $64 + M\log_2(40)$ bits, where the first $64$ bits are used for transmitting $\acute{\lambda}_{k,i} := \max(\abs{\gLki})$, and each $M$ scalars in $\gLki / \acute{\lambda}_{k,i}$ are (dither) quantized to $40$ levels uniformly chosen between $(-1,1]$.
The broadcast is received by any agent in the network with the instantaneous capacity less than the rate.
When an agent $v$ receives $\est{\gL}_{u,i}$ from the broadcasting agent $u$, it updates its estimate via $\vh_{v,i+1} = (1 - \beta_i) \gL_{v,i} + \beta_i \est{\gL}_{u,i}$, where $\beta_i$ is the same as for the OTA-C.
The agents are scheduled to broadcast in a TDMA fashion, i.e.~the next agent starts only after the previous one has finished its broadcast using the $1$MHz bandwidth.
For each broadcast, the receiver is assumed to have perfect CSI knowledge.

In the ANB scheme, similar to the DBC scheme, agents are scheduled to broadcast in a TDMA fashion.
A broadcasting agent $r$ at time $i$ encodes its information $\gL_{r,i}$ into $(M+1)B$ symbols ($MB$ are information bearing symbols and $B$ symbols are for training) using the encoding operation \eqref{eq:encoder}.
Each receiving agent $j\in\sN_r$ receives $(M+1)B$ symbols, and estimates the channel power $\abs{\xi_{jr}}^2$ using the $B$ training symbols.
Then, the agent $j$ obtains the estimate $\est{\gL}_{r,i}$ using remaining $MB$ symbols.
If the received signal power is less than twice the noise power, the received signal is ignored.
Otherwise, the agent $j$ updates its estimate via $\vh_{j,i+1} = (1 - \beta_i) \gL_{j,i} + \beta_i \est{\gL}_{r,i}$, where $\beta_i$ is the same as for the OTA-C.
Therefore, in the low SNR regime, performance of the ANB scheme is at least as good as the NOC scheme.


\subsection{Results and Discussion} \label{sec:sim_results}

\pgfplotsset{every axis/.append style={
                    label style={font=\scriptsize},
                    tick label style={font=\scriptsize},
                    grid=both,
    				grid style=dashed,
    				legend style={nodes={scale=0.7, transform shape}},
                    ylabel near ticks,
                    label style={font=\scriptsize},  
                    }}
\begin{figure}
    \centering
    \begin{tikzpicture}
        \begin{axis}[
        scale only axis=true,
        width=0.85\columnwidth,
        height=0.6\columnwidth,
        xlabel={Iteration $i$},
        ylabel={NMSE (dB)},
        mark size=2.0pt,
        mark repeat=10
        ]
        \addplot [solid, thick, black, mark=*] table[x=Iter, y=NOTAC-U, col sep=comma]{plot_data_nmse.csv};
        \addplot [solid, thick, blue, mark=square*] table[x=Iter, y=NOTAC-W, col sep=comma]{plot_data_nmse.csv};
        \addplot [densely dashed, thick, brown, mark=otimes*] table[x=Iter, y=NBDC-Q40, col sep=comma]{plot_data_nmse.csv};
        \addplot [solid, thick, Rhodamine, mark=star] table[x=Iter, y=APB, col sep=comma]{plot_data_nmse.csv};
        \addplot [densely dotted, thick, red, mark=triangle*] table[x=Iter, y=NNOC, col sep=comma]{plot_data_nmse.csv};
        \addplot [densely dashed, thick, ForestGreen, mark=diamond*] table[x=Iter, y=NCEN, col sep=comma]{plot_data_nmse.csv};
        \legend{OTA-CS, OTA-C, DBC, ANB, NOC, CEN}
        \end{axis}
    \end{tikzpicture}
    \caption{Performance comparison in terms of NMSE for $N=100$.}
    \label{fig:NMSE}
\end{figure}
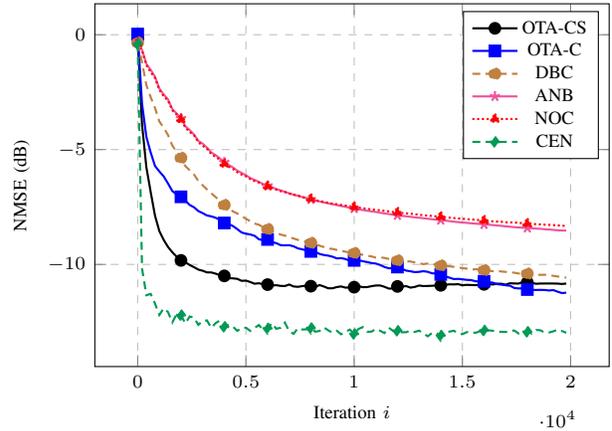
\begin{figure}
    \centering
    \begin{tikzpicture}
        \begin{axis}[
        scale only axis=true,
        width=0.85\columnwidth,
        height=0.6\columnwidth,
        xlabel={Iteration $i$},
        ylabel={NMSE (dB)},
        legend columns=2,
        mark size=1.0pt,
        mark repeat=10
        ]
        \addplot [solid, thick, Gray, mark=*] table[x=Iter, y=NOTAC-N50, col sep=comma]{plot_data_nmse.csv};
        \addplot [dashed, thick, Gray, mark=*] table[x=Iter, y=NBDC-BN50-Q40, col sep=comma]{plot_data_nmse.csv};

        \addplot [solid, thick, Aquamarine, mark=square*] table[x=Iter, y=NOTAC-W, col sep=comma]{plot_data_nmse.csv};
        \addplot [dashed, thick, Aquamarine, mark=square*] table[x=Iter, y=NBDC-Q40, col sep=comma]{plot_data_nmse.csv};

        \addplot [solid, thick, ForestGreen, mark=otimes*] table[x=Iter, y=NOTAC-N150, col sep=comma]{plot_data_nmse.csv};
        \addplot [dashed, thick, ForestGreen, mark=otimes*] table[x=Iter, y=NBDC-BN150-Q40, col sep=comma]{plot_data_nmse.csv};

        \addplot [solid, thick, BrickRed, mark=triangle*] table[x=Iter, y=NOTAC-N200, col sep=comma]{plot_data_nmse.csv};
        \addplot [dashed, thick, BrickRed, mark=triangle*] table[x=Iter, y=NBDC-BN200-Q40, col sep=comma]{plot_data_nmse.csv};

        \addplot [solid, thick, black, mark=diamond*] table[x=Iter, y=NOTAC-N250, col sep=comma]{plot_data_nmse.csv};
        \addplot [dashed, thick, black, mark=diamond*] table[x=Iter, y=NBDC-BN250-Q40, col sep=comma]{plot_data_nmse.csv};

        \legend{OTA-C N=50, DBC N=50, OTA-C N=100, DBC N=100, OTA-C N=150, DBC N=150, OTA-C N=200, DBC N=200, OTA-C N=250, DBC N=250}
        \end{axis}
    \end{tikzpicture}
    \caption{Comparison of OTA-C and DBC schemes in terms of number of agents in the network.}
    \label{fig:nagents}
\end{figure}
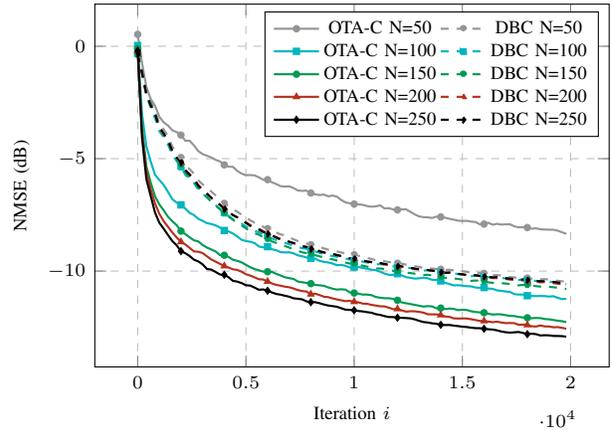
\pgfplotsset{every axis/.append style={
                    label style={font=\scriptsize},
                    tick label style={font=\scriptsize},
                    grid=both,
    				grid style=dashed,
    				legend style={nodes={scale=1.0, transform shape}},
                    ylabel near ticks,
                    label style={font=\scriptsize},  
                    }}
\begin{figure}
    \centering
    \begin{tikzpicture}
        \begin{axis}[
        scale only axis=true,
        width=0.85\columnwidth,
        height=0.6\columnwidth,
        xlabel={SNR (dB)},
        ylabel={NMSE (dB)},
        ymax=-8,
        mark size=2pt,
        legend pos={south west}
        ]
        \addplot [solid, thick, blue, mark=square*] table[x=SNR, y=OTA-C, col sep=comma]{plot_data_snr.csv};
        \addplot [densely dashed, thick, brown, mark=otimes*] table[x=SNR, y=BDC, col sep=comma]{plot_data_snr.csv};
        \addplot [solid, thick, Rhodamine, mark=star] table[x=SNR, y=APB, col sep=comma]{plot_data_snr.csv};
        \addplot [densely dotted, thick, red, mark=triangle*] table[x=SNR, y=NOC, col sep=comma]{plot_data_snr.csv};
        \legend{OTA-C, DBC, ANB, NOC}
        \end{axis}
    \end{tikzpicture}
    \caption{Varying SNR at a distance of $500$m from the transmitter, for $N=100$ at $i=10000$.}
    \label{fig:SNR}
\end{figure}
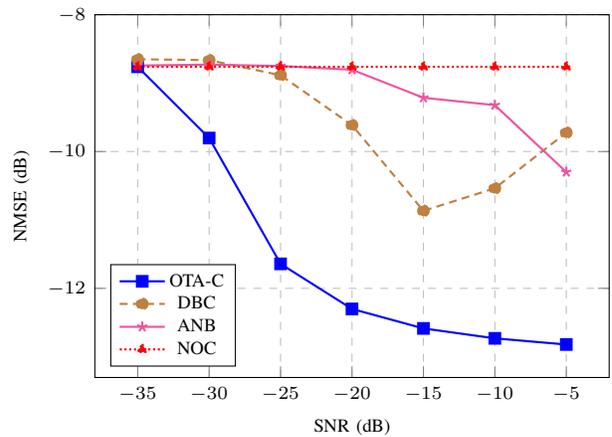

The results are averaged over $100$ independent simulation runs.
For Figures \ref{fig:NMSE} and \ref{fig:nagents}, the SNR is $-20$dB.
Figure \ref{fig:NMSE} shows the performance of different algorithms for $N=100$ agents over $i=0,\dots,20000$ steps in terms of the normalized mean square error (NMSE) in estimating the function $f$ with $\hf$.
The benefits of communication among the agents are evident by comparing the CEN and the NOC schemes.
Here, the ANB performs almost on par with the NOC scheme, since at $-20$dB SNR the noise power is higher than the signal power for most users.
The proposed OTA-C scheme outperforms the separation-based DBC and ANB schemes because the OTA-C scheme enables faster diffusion of information via the consensus protocols.
The sparsity promoting OTA-CS scheme initially converges faster when compared to the vanilla OTA-C scheme, but reaches a crossover point at around $16000$ iterations.
Most notably, the OTA-CS scheme achieves this performance with only $\approx 10\%$ non-zero elements in the parameter vectors.
In other words, even with twice the number of parameters ($2L$ instead of $L$ kernels), the OTA-CS scheme saves up to $80\%$ energy in communication compared to the OTA-C scheme.
Hence, the OTA-CS scheme is suitable for applications that require an energy-efficient method for agents to reach a reasonably good estimate of a solution quickly.
In Figure \ref{fig:nagents}, we compare the performance of the DBC and the OTA-C schemes for different numbers of agents in the network.
As the network size varies from $N=50$ to $N=250$, the performance of the OTA-C scheme improves (almost logarithmically with $N$), while that of the DBC scheme remains nearly the same.

The results in Figure \ref{fig:nagents} demonstrate that the proposed OTA-C scheme easily outperforms the DBC scheme, as the network becomes larger and denser.
To showcase the robustness of the proposed scheme, we compare NMSE performance of the schemes in Figure \ref{fig:SNR} for different SNR values at iteration $i=10000$ for a network with $N=100$ agents.
For the results in Figure \ref{fig:SNR}, the step-size $\beta_i$ is fixed to a small value for low SNR regime, i.e.~$\beta_i := 10^{(\text{SNR} - 15)/20}$ if $\text{SNR}\leq -15$dB, to ensure that at low SNR regimes, the agents fall back to the NOC scheme.
For SNR $>15$dB, we use $\beta_i := 10^{(\text{SNR} - 15)/20} \cdot (\lfloor i/100 \rfloor)^{-0.5}$, to reduce accumulation of the noise with more iterations.
The proposed OTA-C scheme outperforms other schemes at SNRs as low as $-30$dB.
In the DBC scheme, we see a degradation in the performance at SNRs $>-15$dB.
This behavior can be attributed to the diminishing step size $\beta_i$ used in the algorithms for SNRs $\geq -10$dB, which restricts the diffusion of information in the DBC scheme.
On the other hand, in the proposed OTA-C scheme with the same parameters as DBC, a good performance is maintained owing to its superior information diffusion property.
These results also show that the proposed scheme converges when the step-size sequence $(\beta_i)$ is fixed to a reasonably small value.
Overall, the results in Figures \ref{fig:NMSE}, \ref{fig:nagents}, and \ref{fig:SNR} corroborate with our theoretical findings, and provide numerical evidence for efficacy of the proposed scheme.

    \section{Conclusions}
    \label{sec:conc}

In conclusion, this study opens new possibilities for the design and development of communication-efficient and decentralized distributed algorithms. 
The proposed algorithm, with its unique features and compatibility with the superiorization methodology, showcases the potential for achieving high performance in time-varying networks without the need for extensive coordination or knowledge of network specifics.
As we look to the future, there are numerous research avenues that can be explored. 
One direction is the extension of the algorithm to support graphs with random consensus matrices that are not radial (see Definition \ref{def:consensus_matrix}(iv)), offering agents increased autonomy and flexibility, for example, by allowing them to choose their peak transmit powers independently. 
Additionally, leveraging the beamforming capabilities of multiple-input and multiple-output (MIMO) communication systems presents an opportunity to optimize edge weights for even faster convergence. 
Undoubtedly, this work has generated important research questions, and further contributions in this area are expected to advance the field of distributed optimization in wireless networks.

    \bibliographystyle{IEEEtran}
    \bibliography{references}
    
    \clearpage
    
    \section{Supplementary Material} \label{sec:suppl}

The following well-known result is reproduced here, as it is used extensively to prove results of this study.
\begin{proposition}\cite[Theorem 1]{robbins_1971} \cite[Theorem 1]{ermoliev_1983}
    \label{prop:sto_app}
    Let $(\vx_i)_\iN$ be a sequence of random vectors in $\pspace$ such that $\Exp{\norm{\vx_0}^2} < \infty$.
    Suppose that, for a given nonempty set $\sC\subset\Real^M$, any $\vxs \in \sC$, and every $\iN$,
    the following holds a.s.:
    \begin{equation}
        \label{eq:sth_cond_2}
        \begin{aligned}
            \Exp{f(\vx_{i+1}) \big| \vx_i, \dots, \vx_0} \leq (1 + \alpha_i) f(\vx_i) - y_i + z_i,
        \end{aligned}
    \end{equation}
    where we define $f(\vx) := \norm{\vx - \vxs}^2$,
    $(\alpha_i)_\iN$ is a sequence of nonnegative scalars such that $\sif \alpha_i < \infty$,
    and, for all $\iN$, $y_i$ and $z_i$ are nonnegative-valued functions of $\vx_0, \dots, \vx_i$.
    If, in addition, we have $\sif z_i < \infty$,
    then:
    
        (i) the sequence $(\norm{\vx_i - \vxs})_\iN$ converges a.s.~for any $\vxs\in\sC$, and
            $\Exp{\norm{\vx_i - \vxs}^2} < \infty$ for all $\iN$;
            
        (ii) with probability one, $\sif y_i < \infty$;
        
        (iii) the set of accumulation points of $(\vx_i(\omega))_\iN$ is nonempty for a.e.~$\oinO$;
        
        (iv) if two accumulation points $\vx'_i(\omega),\, \vx''_i(\omega)$ of $(\vx_i(\omega))_\iN$
            are such that $\vx'_i(\omega) \notin \sC$ and $\vx''_i(\omega) \notin \sC$ for a.e.~$\oinO$, 
            then the set $\sC$ lies in a hyperplane
            equidistant from $\vx'_i(\omega)$ and $\vx''_i(\omega)$, or, in other words,
            $\norm{\vx'_i(\omega) - \vxs} = \norm{\vx''_i(\omega) - \vxs}$ for all $\vxs\in\sC$.

\end{proposition}

\subsection{Proof of Theorem \ref{thm:main}} \label{apx:thm_main}
To simplify the notation used in \eqref{eq:scheme}, we define $\mLi := \mI - \mPi$, $\mGi := \mI - \beta_i \mLi$, $\mGb := \Exp{\mGi}$, $\mPb := \Exp{\mPi}$, $\mLbi := \Exp{\mLi}$ and $\vti := \mWi \gLi + \vni$.
Using this notation, let us rewrite \eqref{eq:scheme_consensus} as: $(\forall \iN)$ 
\begin{align} \label{eq:protocol}
    \gP_{i+1}   = \mGi \gLi + \beta_i \vti.
\end{align}

As in \cite{cavalcante_dynamic_2013}, we first establish some key properties of the $(\epsilon_0, \delta_0)$-random consensus matrix $\mPi$:
\begin{lemma}
    \label{lem:network_props}
    For all $\iN$, let the matrix $\mPi$ be an $(\epsilon_0, \delta_0)$-random consensus matrix in the sense of Definition \ref{def:consensus_matrix}.
    Then, the following holds for every $\iN$:
    
        (i) $\norm{\mPb}_2 = 1$.
        
        (ii) $\mGi \mJ = \mPi \mJ = \mJ$ and $\mLi\mJ = \vzero$, a.s.
        
        (iii) $\mGb \mJ = \mPb \mJ = \mPb^T \mJ = \mJ$ and $\mLbi\mJ = \mLbi^T \mJ = \vzero$.
        
        (iv) there exists $(\kappa_i)_\iN$ in $\Realp$ such that the series $\sum_{\iN} \kappa_i$ is convergent, and $\norm{\Exp{\mGi^T \mGi}}_2 \leq 1 + \kappa_i$.
            
        (iv) $(\forall \vv\in\VS)$ $\vv^T \mLbi \vv = \vv^T \mLbi^T \vv \geq \epsilon_0 \norm{(\mI - \mJ) \vv}^2$.

    \begin{proof}
        The proofs of parts (ii) and (iii) follow from \cite[Lemma 1]{cavalcante_dynamic_2013}. 
        Parts (i), (iv), and (v) are proved in the following:
        
            (i) By Definition \ref{def:consensus_matrix}(i) and (ii), $\mPb$ is row-stochastic and nonnegative, respectively.
                By Perron's theorem for nonnegative matrices \cite[Theorem 1.2.2]{bapat1997nonnegative}, the spectral radius of $\mPb$ is equal to $1$.
                Since $\mPb$ is a radial matrix, by Definition \ref{def:consensus_matrix}(iv), we have $\norm{\mPb}_2 = \rho(\mPb) = 1$.
        
        
            
            (iv) Using $\mGi = (\mI - \beta_i(\mI - \mPi))$ and $\norm{\mPb}_2 = 1$ (Lemma \ref{lem:network_props}(i)), we obtain the bound $\norm{\Exp{\mGi^T \mGi}}_2 \leq 1 - \beta_i^2 + \beta_i^2 \norm{\Exp{\mPi^T\mPi}}_2$.
                Defining $\kappa_i = \beta_i^2 \delta_0$, the result follows from Definition \ref{def:consensus_matrix}(v).
                The series $\sif \kappa_i = \delta_0 \sif \beta_i^2$ is convergent by definition of $(\beta_i)_\iN$ [see Assumption \ref{ass:local_scheme}(iv)].
                
            (v) 
                Note that $(\mI - \mJ) (\mI - \mPb) (\mI - \mJ) = \mI - \mPb = \mLbi$ by Lemma \ref{lem:network_props}(iii).
                For any $\vv\in\VS$:
                \begin{equation*}
                    \begin{aligned}
                        \vv^T \mLbi \vv 
                            &= \vv^T (\mI - \mJ) (\mI - \mPb) (\mI - \mJ) \vv \\
                            &= \vv^T (\mI - \mJ) (\mI - \mJ) \vv - \vv^T (\mI - \mJ) \mPb (\mI - \mJ) \vv \\
                            & \mymark{(a)}{\geq} \norm{(\mI - \mJ) \vv}^2 - (1 - \epsilon_0) \norm{(\mI - \mJ) \vv}^2 \\
                            &= \epsilon_0 \norm{(\mI - \mJ) \vv}^2,
                    \end{aligned}
                \end{equation*}
                where $(a)$ follows from Definition \ref{def:consensus_matrix}(iii), and the fact that for any $\vv\in\VS$, $(\mI - \mJ) \vv \in \sC^{\perp}$.
    \end{proof}
\end{lemma}

        


\subsubsection{Proof of Theorem \ref{thm:main}(i)}

    The proof is based on the application of Proposition \ref{prop:sto_app} to the sequence $(\norm{\gP - \gPs})_\iN$ by establishing the sufficient conditions, in particular, condition \eqref{eq:sth_cond_2}.
    We define two dummy variables for the following discussions $\vy_i := \opTi(\gPi - \gHi + \zeta_i\gZi) - \gPs$, and $\fs(\gP) := \norm{\gP - \gPs}^2$, where $\gPs:=[(\bh)^T, \dots, (\bh)^T]^T$ for some $\bh\in\sUs$.
    Using \eqref{eq:protocol}, we can expand $\fsp$ to get:
    $\fsp = \norm{\gP_{i+1} - \gPs}^2 = \norm{\mGi \vy_i + \beta_i \vti}^2 = \vy_i^T \mGi^T \mGi \vy_i + \beta_i^2 \norm{\vti}^2 + 2 \beta_i \vy_i^T\mGi^T \vti$.
    Now, consider the conditional expectation $\Epp{\fsp}$:
    \begin{align}
        &\Epp{\fsp} \mymark{(a)}{=} \Epp{\vy_i^T \Exp{\mGi^T \mGi} \vy_i} \nonumber\\
        &\quad\qquad + \beta_i^2\ \Epp{\norm{\vti}^2} + 2 \beta_i \Epp{\vy_i^T \mGi^T \vti} \nonumber\\
        &\ \mymark{(b)}{\leq} \norm{\Exp{\mGi^T \mGi}}_2 \Epp{\norm{\vy_i}^2} + \beta_i^2\ \Epp{\norm{\vti}^2} \nonumber\\
        &\  \qquad+ 2 \beta_i \Epp{\vy_i^T \mGi^T \vti} \nonumber\\
        &\ \mymark{(c)}{\leq} (1 + \kappa_i)\ \Epp{\norm{\vy_i}^2} + \beta_i^2\ \Epp{\norm{\vti}^2}. \label{eq:apx_i_2}
    \end{align}
    In (a), we use the law of total expectation, whereby, for any random variable (r.v.) $x$ that depends on the r.v.~$\gPi$ and $\gZi$, we have $\Epp{x} = \Expp{\gZi}{\Exp{x \big| \gPi, \gZi} \big| \gPi}$.
    The step (b) follows from the spectral norm inequality.
    The first term in (c) follows from Lemma \ref{lem:network_props}(iv).
    We show that $\Epp{\vy_i^T \mGi^T \vti} = 0$ by expanding it as follows:
    \begin{align*}
        &\Epp{\vy_i^T \mGi^T \vti} = \Epp{(\gLi - \gPs)^T \mGi^T (\mWi \gLi + \vni)} \\
        &\ = \mathbb{E}\Big[\gLi^T \Exp{\mGi}^T \Exp{\mWi} \gLi + \gLi^T \Exp{\mGi}^T \Exp{\vni} \\
        &\qquad  - (\gPs)^T \Exp{\mGi}^T \Exp{\mWi} \gLi + \Exp{\mGi}^T \Exp{\vni} \big| \gPi \Big].
    \end{align*}
    Each summand inside the outer expectation is zero due to the following facts: (1) $\mGi$ is uncorrelated to $\mWi$ and $\vni$, and (2) $\mWi$ and $\vni$ are zero mean [see Assumptions \ref{ass:random_ass}(ii) and \ref{ass:random_ass}(iii)].

    From Assumption \ref{ass:bounded_subg}, and the fact that $\opLi$ in \eqref{eq:comm_model} is a bounded linear operator, we verify that the sequence $(\gPi - \gHi - \gPs)_\iN$ is bounded a.s.~(recall that $\gPs$ is fixed).
    Hence, for a.e.~$\oinO$, there exists $u_1(\omega)\in\Realp$ such that for all $\iN$, we have $\norm{\gPi - \gHi - \gPs} \leq u_1 := u_1(\omega)$.
    
    In the following, we upper-bound the terms $\Epp{\norm{\vy_i}^2}$ and $\Epp{\norm{\vti}^2}$ in \eqref{eq:apx_i_2}:
    
        (A) \ul{Bounding $\Epp{\norm{\vy_i}^2}$}:
        Since $\opTi$ is a quasi-nonexpansive operator [Assumption \ref{ass:local_scheme}(i)] and $\gPs$ is an element of $\Fix{\opTi}$, we have $\Epp{\norm{\vy_i}^2} \leq \Epp{\norm{\gPi - \gHi + \zeta_i\gZi - \gPs}^2}$.
        Expanding the right-hand side of the inequality, we get:
        \begin{align}
            &\Epp{\norm{\vy_i}^2} \leq  \fse + \norm{\gHi}^2 - 2(\gHi)^T\ (\gPi - \gPs) \nonumber\\
                &\quad\qquad    + \zeta_i^2\Epp{\norm{\gZi}^2} + 2\zeta_i\Epp{\gZi^T}(\gPi - \gHi - \gPs) \nonumber\\
            &\quad \mymark{(a)}{\leq} \fse - \sum_{k\in\sN} \rmuki (2 - \rmuki) \frac{(\cT(\vhki))^2}{\norm{\cS(\vhki)}^2} \nonumber\\
                &\quad\qquad + \zeta_i^2 N r^2 + 2 \zeta_i \sqrt{N} r u_1.        \label{eq:apx_i_6}
        \end{align}
        In $(a)$ above, the second term on the r.h.s.~follows by using the definition of subgradient in \eqref{eq:alpha} (for details see \cite[Theorem 2(a)]{yamada2005adaptive}).
        The last two terms of the r.h.s.~in $(a)$ is a consequence of $(\zeta_i \vzki)_\iN$, for all $\kN$ and $\iN$, being a sequence of bounded perturbations in the sense of Definition \ref{def:perturb}. 
        More precisely, for every $\iN$, there exists a scalar $r\in\Realp$, independent of $i$, such that $\Epp{\norm{\gZi}^2} = \sum_{\kN} \Epp{\norm{\vzki}^2} \leq N r^2$.
        
        (B) \ul{Bounding $\Epp{\norm{\vti}^2}$}: 
        Using the law of total expectation, and the fact that, given r.v.~$\gPi$ and $\vz_i$, the r.v.~$\gLi$ can be completely determined, we have:
        \begin{align}
            &\Epp{\norm{\vti}^2} = \Epp{\norm{\mWi \gLi + \vni}^2} \nonumber\\
            &\ \leq \Epp{\norm{\gLi}^2} \norm{\Exp{\mWi^T \mWi}}_2 + \Exp{\norm{\vni}^2}  \nonumber\\
            &\ \qquad + 2\Epp{\norm{\gLi}} \norm{\Exp{\mWi^T \vni}} \nonumber\\
            &\ \leq a_0 u_2^2 + b_0 + 2 c_0 u_2 =: \varrho < \infty, \label{eq:apx_1_4}
        \end{align}
        where, in the last inequality, scalars $a_0, b_0,$ and $c_0$ are defined in Assumption \ref{ass:random_ass}(iv), and it follows from Assumption \ref{ass:bounded_subg} that there exists $u_2\in\Realp$ such that $\Epp{\norm{\gLi}^2} \leq u_2$ for every $\iN$.

    
    
    \vspace{0.5em}
    Using \eqref{eq:apx_i_6} and \eqref{eq:apx_1_4}, equation \eqref{eq:apx_i_2} can be written as:
        \begin{align} \label{eq:apx_i_7}
            \Epp{\fsp} \leq (1 + \kappa_i) \fse + \varphi_i - \vartheta_i,
        \end{align}
    where 
    $\varphi_i := \left(\varrho\beta_i^2 + (1 + \kappa_i)(\zeta_i^2 Nr^2 + 2 \zeta_i \sqrt{N}ru) \right)$, and
    $\vartheta_i := (1 + \kappa_i) \rho_i$, where $\rho_i := \sum_{k\in\sN} \rmuki (2 - \rmuki) (\cT(\vhki))^2 / \norm{\cS(\vhki)}^2$.
    Note that the following: (a) series $\sum_{\iN} \kappa_i$ is convergent, which follows from Lemma \ref{lem:network_props}(iv);
    (b) series $\sum_\iN \varphi_i$ is convergent because $\sum_\iN \beta_i^2$, $\sum_\iN \kappa_i$ and $\sum_\iN \zeta_i$ are convergent;%
    \footnote{Here, we use the fact that for any two convergent series $\sum_\iN a_i$ and $\sum_\iN b_i$ of nonnegative scalars, $\sum_\iN a_i b_i$ is convergent.}~%
    and (c) $\vartheta_i$ is nonnegative due to Assumption \ref{ass:local_scheme}(ii).

    Hence, the conditions are satisfied for application of Proposition \ref{prop:sto_app} to \eqref{eq:apx_i_7}, from which we can conclude that the following results hold a.s.:
    
    (a) The sequence $\{\norm{\gPi - \gPs}^2\}$ converges.
    
    (b) The series $(\vartheta_i)_\iN$ converges, which, based on $\rmuki \in (0,2)$, and that series $(\kappa_i)$ is convergent, implies that the series
            \begin{align}
                \label{eq:subg_ratio_inq}
                \sum_{\iN} \sum_{k\in\sN} 
                    \frac{(\cT(\vhki))^2}{\norm{\cS(\vhki)}^2}
            \end{align}
            converges.
            In particular, since the sequence $(\cS(\vhki))_\iN$ is bounded by Assumption \ref{ass:bounded_subg}, so convergence of \eqref{eq:subg_ratio_inq} also shows that $\lim_{i\to\infty} \cT(\vhki) = 0$ for all $k\in\sN$.
    This concludes the proof of Theorem \ref{thm:main}.
    
\vspace{1em}

\subsubsection{Proof of Theorem \ref{thm:main}(ii)}

    Using $\mGi := \mI - \beta_i \mLi$, expand the expression for $\Epp{\fsp} = \Epp{\norm{\gPi - \gPs}^2}$ to get:
%
    \begin{align}
        &\Epp{\fsp} = \Epp{\norm{\vy_i}^2} - 2\beta_i \Exp{\vy_i^T \mLi \vy_i} \nonumber\\
        &\quad\quad + \beta_i^2 \Epp{\norm{\vti}^2 \Big| \ \gPi} + \beta_i^2 \Epp{\vy_i^T \Exp{\mLi^T\mLi} \vy_i} \nonumber\\
        &\quad\leq \Epp{\norm{\vy_i}^2} - 2\beta_i \Epp{\vy_i^T \mLi \vy_i} + \nonumber\\
        &\quad  \beta_i^2 \norm{\Exp{\mLi^T\mLi}}_2 \Epp{\norm{\vy_i}^2} + \beta_i^2 \Epp{\norm{\vti}^2}. \label{eq:apx_ii_x}
    \end{align}
    
    The following inequalities hold for terms in \eqref{eq:apx_ii_x}, for all $\iN$:
    
    (1) $\Epp{\norm{\vy_i}^2} \leq \fse + \zeta_i^2 Nr^2 + \zeta_i \sqrt{N}ru$, which follows from \eqref{eq:apx_i_6}, and because $\rho_i \geq 0$ [see below \eqref{eq:apx_i_7}].
    
    (2) $\norm{\Exp{\mLi^T \mLi}}_2 \leq 4 + \delta_0$, where, using $\mLi = \mI - \mPi$ in the l.h.s., it follows from Lemma \ref{lem:network_props}(i) (i.e.~$\norm{\mPb}_2=1$), and $\norm{\Exp{\mPi^T\mPi}}_2 \leq 1 + \delta_0$ from Definition \ref{def:consensus_matrix}(v).
            
    (3) $\Epp{\vy_i^T \mLi \vy_i} \geq \epsilon_0\ \Epp{\norm{(\mI - \mJ) \gLi}^2}$ follows from the law of total expectation, Lemma \ref{lem:network_props}(v), and the fact that $(\mI - \mJ) \gPs = \mZero$. 

    (4) $\Exp{\norm{\vti}^2} \leq \varrho$, where $\varrho\in\Realp$ is defined in \eqref{eq:apx_1_4}.


    Using the inequalities above, we can rewrite \eqref{eq:apx_ii_x} as:
    \begin{align}
        & \Epp{\fsp} \leq (1 + \check{\kappa}_i)\fse + \check{\varphi}_i - \check{\vartheta}_i, \label{eq:apx_ii_2}
    \end{align}
    where we define $\check{\kappa}_i := \beta_i^2 (4 + \delta_0)$, $\check{\varphi}_i := (\beta_i^2 \big(\zeta_i^2 Nr^2 + \zeta_i \sqrt{N}ru + \varrho \big) + \zeta_i^2 Nr^2 + \zeta_i \sqrt{N}ru)$, and $\check{\vartheta}_i := 2 \epsilon_0 \beta_i \Epp{\norm{(\mI - \mJ) \gLi}^2}$.
    Both series $\check{\kappa}_i$ and $\check{\varphi}_i$ are convergent by definition.
    Hence, the conditions are satisfied for application of Proposition \ref{prop:sto_app} to \eqref{eq:apx_ii_2}.
    In addition to the result already established in Theorem \ref{thm:main}(i), 
    we verify that the series $(\check{\vartheta}_i)_\iN$ is convergent.
    Since $\sum_{\iN} \beta_i$ diverges to infinity, and $\epsilon_0 > 0$ by Assumption \ref{ass:random_ass}(i), we obtain: 
    \begin{align} \label{eq:apx_ii_y}
        \liminf\limits_{i\to\infty} \Epp{\norm{(\mI - \mJ) \gLi}^2} = 0 \quad \text{ a.s.}
    \end{align}
    Equation \eqref{eq:apx_ii_y} implies that there is a subsequence $l\in\sK$, $\sK\subset \Natural$, such that the sequence of random variables $(\Epp{\norm{(\mI - \mJ) \gL_l}^2})_\iN$ converge to zero a.s. 
    Now, using \eqref{eq:scheme_consensus}, we obtain the inequality: $\Epl{\norm{(\mI - \mJ) \gP_{l+1}}^2}\leq (1 - \beta_l) \Epl{\norm{(\mI - \mJ)\gL_l}} + \beta_l \Epl{\norm{ (\mI - \mJ) \opL_l(\gL_l)}}$.
    Hence, by using the facts that (1) $\limsup_{l\to\infty} \beta_l = 0$, and (2) the random linear operator $\opLi$ maps any bounded vector $\gLi$ to a bounded set, it can be verified that the limit $\lim_{l\to\infty} \Epl{\norm{(\mI - \mJ) \gP_{l+1}}} = 0$ exists.
    This proves Theorem \ref{thm:main}(ii).

\subsubsection{Proofs of Theorem \ref{thm:main}(iii) and \ref{thm:main}(iv)}
    Equation \eqref{eq:apx_ii_y} above, along with Assumption (b) in Theorem \ref{thm:main}(iii), i.e.~$\Epp{\gLi} = \gLi$ a.s., verifies that $\liminf\limits_{i\to\infty} \norm{(\mI - \mJ) \gLi} = 0$ a.s.
    With this result in place, the proofs of Theorem \ref{thm:main}(iii) and \ref{thm:main}(iv) follows from the proofs of similar statements in \cite[Theorem 2]{cavalcante_dynamic_2013}.

\subsection{Proof of Lemma \ref{lem:unbiased_est}} \label{apx:lem_otac}
By expanding the terms $\norm{\vy_r}^2$ and $\norm{\vy'_r}^2$ in \eqref{eq:r_postproc} using \eqref{eq:WMAC}, it follows from simple algebraic manipulations that
\begin{align*}
    \est{v}_{r} &= \sum_{j\in\sN_r} \left((P_j/B) \sum_{n=1}^B \abs{\xi_{jr}(2n-1)}^2\right) \lambda_j + \eta_r, \\
    \est{v}'_r &= \sum_{j\in\sN_r} \underbrace{\left((P_j/B) \sum_{n=1}^B \abs{\xi_{jr}(2n-1)}^2\right)}_{=:\ \rnu_{jr}}\ +\ \eta'_r,
\end{align*}
where $\eta_r$ and $\eta'_r$ are given by:
\begin{equation}
    \label{eq:OTA_noise}
    \begin{aligned}
        \eta_r  = \Delta (\tW_r - \eNoise) + \dmin \eta'_r,  \quad
        \eta'_r = \tW'_r - \eNoise,
    \end{aligned}
\end{equation}
where $\tW_r$ is given by:
\begin{align*}
    \tW_r &:= \frac{1}{B} \left(\sum_{b=1,3,\dots}^{2B-1} \tau_r(b) + \abs{w_r(b)}^2\right),
\end{align*}
\begin{align*}
    \tau_r(b) &:= \sum_{\substack{k,l \\ k \neq l}} \sqrt{g_k(\lambda_k) g_l(\lambda_l)}\ \xi_{kr}(b) \xi_{lr}^*(b)\ U_k(b) U_l(b) \\
        &+ 2 \RE{w_r^*(b) \sum_{k\in\sN_r} \sqrt{g_k(\lambda_k)}\ \xi_{kr}(b)\ U_k(b)}
\end{align*}
and $\tW'_r$ is defined similarly for $b=2, 4, \dots, 2B$ with $\lambda_k,\lambda_l$ replaced by $\dmax$ for all $k,l\in\sN_r$.
It follows that $\Exp{\eta_r} = \Exp{\eta'_r} = 0$ because (i) $(U_k(b))$ are zero-mean and i.i.d.~random variables (r.v.), (ii) the noise r.v.~$w_r$ is zero-mean and independent of the channel [Assumption \ref{ass:WMAC}(i) and (vi)], and (iii) $\Exp{\tW_r} = \Exp{\tW'_r} = \eNoise := \Exp{\abs{w_r}^2}$.

\subsection{Proof of Proposition \ref{prop:OTA_scheme}}
\label{apx:prop_ota}


The proof essentially follows from Assumption \ref{ass:WMAC} on the WMAC model, and Assumption \ref{ass:network} on the network graph.
First, we prove the condition in Assumption \ref{ass:random_ass}(i), i.e., for all $\iN$, there exist $0<\epsilon_0\leq 1$ and $0\leq \delta_0 < \infty$ such that $\oPi$ is an $(\epsilon_0, \delta_0)$-random consensus matrix in the sense of Definition \ref{def:consensus_matrix}.
In this direction, note that the conditions (i)--(iv) in Definition \ref{def:consensus_matrix} follow directly from the definition of $\oPi$ in \eqref{eq:matrix_P} and Assumption \ref{ass:network}.
In particular, condition (iii) in Definition \ref{def:consensus_matrix} follows from the connectivity of the expected graph $\sGb$ (Assumption \ref{ass:network}(i)) \cite{chai_directed_graph_2007}, and condition (iv) follows since $\oPi$ is symmetric because the network graph is undirected (Assumption \ref{ass:network}(ii)).
In the following, we prove the condition (v) in Definition \ref{def:consensus_matrix}.
Using Jensen's inequality, it follows that $\norm{\Exp{\oPi^T\oPi}}_2 \leq \Exp{\norm{\oPi}_2^2}$.
For all $k\in\sN$ and $m=1,\dots,M$, define r.v.~$\tZ_{(k-1)M+m, i} := \sum_{j\in\sN_k} \rnu_{jk,i}^{(m)}$, and $Y_i := \max\{\tZ_{1, i}, \dots, \tZ_{NM, i}\}$.
In words, r.v.~$Y_i$ selects the maximum from the diagonal of matrix $\gamma_i^{-1}(\mI_{MN} - \oPi)$.
If the r.v.~$Y_i \leq \gamma_i^{-1}$ [defined in \eqref{eq:gamma}], then the random matrix $\oPi$ is nonnegative by definition \eqref{eq:matrix_P}.
In this case, since $\oPi$ is row-stochastic, by Gershgorin's circle theorem \cite[Corollary 6.1.5]{horn2012matrix}, the spectral radius of $\oPi$ is equal to $1$.
Otherwise, if $Y_i > \gamma_i^{-1}$, random matrix $\oPi$ may have negative elements on its diagonal, and, again by Gershgorin's circle theorem, the spectral radius is bounded above by the largest sum of the absolute values of its rows.
More precisely, for a given outcome $y_i \sim Y_i$, we have
\begin{equation*}
    \begin{aligned}
        \norm{\oPi}_2 = \rho(\oPi) \leq \begin{cases}
            1, & \text{ if } y_i \leq \frac{1}{\gamma_i} \\
            2\gamma_i y_i - 1, & \text{ otherwise}.
            \end{cases}
    \end{aligned}
\end{equation*}
Let $f$ be the pdf of the r.v.~$Y_i$. 
Then, we can bound the expected value of $\norm{\oPi}_2^2$ as follows:
\begin{equation*}
    \begin{aligned}
        &\Exp{\norm{\oPi}_2^2} 
            \leq \int_{0}^{\frac{1}{\gamma_i}} f(y_i) \deriv y_i + \int_{\frac{1}{\gamma_i}}^\infty (2\gamma_i y_i - 1)^2 f(y_i) \deriv y_i \\
            &\quad \leq 1 + 4 \gamma_i^2 \int_{\frac{1}{\gamma_i}}^\infty y_i^2 f(y_i) \deriv y_i \leq 1 + 4 \gamma_i^2 \Expp{Y_i}{y_i^2}. 
    \end{aligned}
\end{equation*}
Here, $4 \gamma_i^2 \Expp{Y_i}{y_i^2} =: \delta_0$ is positive and finite, since both $\gamma_i$ and $\Expp{Y_i}{y_i^2}$, the second moment of $Y_i$, are positive and finite [see Assumption \ref{ass:WMAC}].

The conditions in Assumption \ref{ass:random_ass}(ii)--(iv) can be verified using definitions of $\oPi$, $\oWi$, and $\ovni$, and Assumption \ref{ass:WMAC}.
This concludes the proof of Proposition \ref{prop:OTA_scheme}.

\subsection{Proof of Proposition \ref{prop:perturb}} \label{apx:prop_sim}
Since $(\forall \iN)\ \zeta_i \in \Realp$ and the series $\sum_\iN \zeta_i$ is convergent by definition, it only remains to prove that the random sequence $(\vzki)_\iN$ is bounded a.s.~for all $\kN$.
With this aim, notice that by replacing $\fq_{k,i}^{(m)}$ into \eqref{eq:perturb_seq}, the $m$th component of $\vzki$, for $m=1,\dots,M$, can be bounded as $\abs{\vzki[m]} \leq (\abs{\vy_{k,i}[m]} + \varsigma)^{-1} \leq \varsigma^{-1}$.
As $\varsigma$ is strictly positive, each $\abs{\vzki[m]}$ for $m=1,\dots,M$, is bounded above.
Hence, $\norm{\vzki}^2 = \sum_{m=1}^M \abs{\vzki[m]}^2$ is bounded above for all $\kN$ and $\iN$.
This completes the proof of the statement in Proposition \ref{prop:perturb}.

\end{document}